\documentclass[11pt]{article} 
\usepackage{latexsym}
\usepackage{amsmath, rotating}
\usepackage{a4wide}
\usepackage{amsmath, rotating, color}
\usepackage{amsfonts,amssymb, amsthm}
\usepackage{psfrag}
\usepackage{pictex}
\newcommand\unnumberedfootnote[1]{ %
        \let\temp=\thefootnote %
        \renewcommand{\thefootnote}{}%
        \footnote{#1}%
        \let\thefootnote=\temp%
        \addtocounter{footnote}{-1}}

\newcommand{\eps}  {\varepsilon}

\newcommand{\R}{\mathbb{R}}
\renewcommand{\P}{\mathbb{P}}
\newcommand{\E}{\mathbb{E}}

\def \bpf {\noindent{\sc Proof: }}
\def \epf {\hbox{}\nobreak\hfill
\vrule width 2mm height 2mm depth 0mm
\par \goodbreak \smallskip}

\newtheorem{theorem}{Theorem}
\newtheorem{proposition}{Proposition}[section]
\newtheorem{lemma}[proposition]{Lemma}
\newtheorem{corollary}[proposition]{Corollary}
\newtheorem{definition}[proposition]{Definition}
\theoremstyle{definition}
\newtheorem{remark}[proposition]{Remark}

\numberwithin{equation}{section}
\date{}
\begin{document}
    \title{\LARGE ``Trees under attack'': a Ray-Knight representation of \\ Feller's branching  diffusion with logistic growth} \author{{\sc by
    V. Le, E. Pardoux and A. Wakolbinger} \\[2ex]
  \emph{Universit\'e de Provence} \\
  \emph{and Goethe--Universit\"at Frankfurt}}
\maketitle
\unnumberedfootnote{\emph{AMS 2000 subject classification.} {\tt 60J70}
 (Primary) {\tt 60J80, 60J55} (Secondary).}
\unnumberedfootnote{\emph{Keywords and phrases.}  {\tt Ray-Knight representation, Feller branching with logistic growth, exploration process, local time, Girsanov transform.}}

\begin{abstract}
We obtain a representation of Feller's branching diffusion with logistic growth in terms of the local times of a reflected Brownian motion $H$  with a drift that is affine linear in the local time accumulated by $H$ at its current level. As in the classical Ray-Knight representation, the excursions of $H$ are the exploration paths of the trees of descendants of the ancestors at time $t=0$, and the local time of $H$ at height $t$ measures the population size at time $t$ (see e.g. \cite{LG4}). We cope with the dependence in the reproduction by introducing a pecking order of individuals: an individual explored at time $s$ and living at time $t=H_s$ is prone to be killed by any of its contemporaneans that have been explored so far. The proof of our main result relies on approximating $H$ with a sequence of Harris paths $H^N$ which figure in a Ray-Knight representation of the total mass of a branching particle system. We obtain a suitable joint convergence of $H^N$ together with its local times {\em and} with the Girsanov densities that introduce the dependence in the reproduction.
\end{abstract}

\section{Introduction}
{\em Feller's branching diffusion with logistic growth} is  is governed by the SDE
\begin{align}\label{fellog1}
dZ_t= \sigma \sqrt{Z_t}\, dW_t \, +(\theta Z_t -\gamma Z_t^2)\, dt,  \qquad   Z_0 =x > 0,
\end{align}
with positive constants $\sigma$, $\theta$ and $\gamma$. It has been studied in detail by Lambert \cite{AL}, and models the evolution of the size of a large population with competition. The diffusion term in \eqref{fellog1} incorporates the individual offspring variance, and the drift term includes a super--criticality in the branching that is counteracted by a killing with a rate proportional to the ``number of pairs of individuals''.

For $\theta = \gamma =0$, equation \eqref{fellog1} is the SDE of {\em Feller's critical branching diffusion with variance parameter $\sigma^2$}. In this case, a celebrated theorem due to Ray and Knight says that $Z$ has a representation in terms of the local times of reflected Brownian motion. To be specific, let $H = (H_s)_{s\ge 0}$ be a Brownian motion on $\mathbb R_+$ with variance parameter $4/\sigma^2$, reflected at the origin, and for $s,t \ge 0$ let $L_s(t,H)$ be the (semi--martingale) local time accumulated by $H$ at level $t$ up to time $s$. Define 
\begin{align}\label{Sx}
S_x := \inf\{s>0: (\sigma^2/4) L_s(0,H) > x\}.
\end{align}
 Then $(\sigma^2/4)  L_{S_x}(t,H), \, t\ge 0,$ is a weak solution of \eqref{fellog1} with $\theta = \gamma =0$ , and is called the {\em Ray-Knight representation} of Feller's critical branching diffusion.

The Ray-Knight representation has a beautiful interpretation in an individual-based picture. Reflected Brownian motion $H = (H_s)_{s\ge 0}$ arises as a concatenation of excursions, and each of these excursions codes a {\em continuum random tree}, the genealogical tree of the progeny of an individual that was present at time $t=0$. The size of this progeny at time $t>0$ is $\sigma^2/4$ times the (total) local time spent by this excursion at level $t$.
Starting with mass $x$ at time $t=0$ amounts  to collecting a local time $(4/\sigma^2) x$ of $H$ at level $0$. The local time of $H$ at level $t$ then arises as a sum over the local time of the excursions, just as the state at time $t$  of Feller's branching diffusion, $Z_t$, arises as a sum of the masses of countably many families, each of which belongs to the progeny of one ancestor that lived at time $t=0$. The path $(H_s)_{0\le s\le S_x}$ can be viewed as the {\em exploration path of the genealogical forest} arising from the ancestral mass $x$. We will briefly illustrate this in Section \ref{discretemass}
along a discrete mass -- continuous time approximation. For a more detailed explanation and some historical background we refer to the survey \cite{PW}.

The motivation of the present paper was the question whether a similar picture is true also for \eqref{fellog1} with strictly positive $\theta$ and $\gamma$, and whether also in this case a Ray-Knight representation is available for a suitably re-defined dynamics of an exploration process  $H$. At first sight this seems prohibiting since the nonlinear term in \eqref{fellog1} destroys the independence in the reproduction. However, it turns out that introducing an order among the individuals helps to overcome this hurdle. We will think of the individuals as being arranged ``from left to right'', and decree that the pairwise fights are always won by the individual ``to the left'', and lethal for the individual ``to the right''. 
In this way we arrive at a population dynamics which leaves the evolution \eqref{fellog1} of the total mass unchanged, see again the explanation in Section \ref{discretemass}. The death rate coming from the pairwise fights leads in the exploration process of the genealogical forest to a downward drift which is proportional to $L_s(H_s, H)$, that is, proportional to the amount of mass seen to the left of the individual encountered at exploration time $s$ (and living at real time $H_s$) -- more rigorously, $L_r(t,H)$ denotes the local time accumulated by the semimartingale $H$ up to time $r$ at level $t$, see the beginning of Section \ref{sec2} for a precise definition.
As a consequence, those excursions of $H$ which come later in the exploration tend to be smaller - the trees to the right are ``under attack from those to the left''.

In quantitative terms, we will consider the stochastic differential equation  
\begin{equation}\label{reflecSDE}
      H_s=\frac{2}{\sigma}B_s+\frac{1}{2}L_s(0,H)+\frac{2\theta }{\sigma^2}s- \gamma \int_0^sL_r(H_r,H)dr, \qquad s \ge 0,
      \end{equation}
where $B$ is a standard Brownian motion. 
The last two terms are the above described components of the drift in the exploration process, and the term $L_s(0,H)/2$ takes care of the reflection of $H$ at the origin.
We will show
\begin{proposition}\label{exun}
The SDE \eqref{reflecSDE} has a unique weak solution.
\end{proposition}
Our main result is the 
\begin{theorem}\label{main}
Assume that $H$ solves the SDE \eqref{reflecSDE}, and
let, for $x>0$, $S_x$ be defined as in \eqref{Sx}. Then  $(\sigma^2/4)L_{S_x}(t,H)$, $t\ge 0$, solves  \eqref{fellog1}.
\end{theorem}
  We will prove Proposition \ref{exun} by a Girsanov argument, and Theorem \ref{main} along a discrete mass--continuous time approximation that is presented in Section \ref{discretemass}.  In section \ref{secZ} we take the limit in the total population process along the discrete mass approximation. An important step in the proof of Theorem \ref{main}, and interesting in its own right, is Theorem \ref{conHL} in Section \ref{allzero}, in which we obtain a convergence in distribution, in the case $\theta=\gamma=0$,  of processes that approximate reflected Brownian motion, together with their local times (considered as random fields indexed by their two parameters, time and level). A similar 
  convergence result was proved in \cite{Pe} for piecewise linear interpolation of discrete time random walks and their local time. The proof of Theorem \ref{main} is completed in Section \ref{secgirs}, using again Girsanov's theorem, this time for Poisson point processes. In that section we will also prove Theorem \ref{th3}, which says that in the case $\theta,\gamma>0$, the exploration process {\em together} with its local time (now evaluated at a certain random time, while the parameter for the various levels is varying) converges along the discrete mass approximation. Finally, for the convenience of the reader, we collect in an Appendix, at the end of this paper, several results from the literature, concerning tightness and weak convergence in the
  space $D$, the Dol\'eans exponential and ``goodness'', and the two versions of Girsanov's theorem which we need : the one for Brownian motion and the one for Poisson point processes.
  
 When this work was already completed, our attention was drawn by Jean--Fran\c cois Le Gall on the article \cite{NRW} by J. Norris, L.C.G. Rogers and D. Willams, who proved a Ray--Knight theorem  for a Brownian motion with a ``local time drift'', using tools from stochastic analysis, in particular the ``excursion filtration''. With a similar methodology we were recently able to establish another, shorter but less intuitive, proof of the main result of this paper \cite{PW2}.

 
 \section{Proof of Proposition \ref{exun}} \label{sec2}
 For abbreviation we write $L_s(t) := L_s(t,H)$ for the local time at level $t\ge0$, accumulated by the continuous semi--martingale 
 $H$ up to time $s$. For the reader's convenience, we recall a possible definition (borrowed from \cite{ReYo}, Theorem VI.1.2)
 of that quantity~:
 $$L_s(t):=2(H_s-t)^+-2\int_0^s{\bf1}_{\{H_r>t\}}dH_r.$$

Following the Girsanov route, let $B$ be a standard Brownian motion defined on a probability space $(\Omega, \mathcal F, \mathbb P)$ and let $H$ obey 
 \begin{equation}\label{BR}
H_s=\frac{2}{\sigma}B_s+\frac{1}{2}L_s(0),
\end{equation}
which can be viewed as Skorohod's equation for reflected Brownian motion, see section 3.6.C in \cite{KS}.
Thus the process $H$ is under the probability measure $\mathbb P$ a reflected Brownian motion with variance parameter $4/\sigma^2$. In order to prove the existence part of Proposition \ref{exun} it suffices to construct, by a suitable reweighting of $\mathbb P$, a probability measure $\tilde  {\mathbb P}$ under which 
\begin{equation}\label{Btilde}
\tilde {B}_s:=B_s-\int_0^s\left[\frac{\theta}{\sigma}-\frac{\sigma\gamma}{2}L_r(H_r)\right]dr, \qquad s\ge 0,
\end{equation}
 is a standard Brownian motion. Indeed, \eqref{Btilde} and \eqref{BR} together imply that under $\tilde {\mathbb P}$ the process $H$ solves \eqref{reflecSDE} with $\tilde B$ instead of $B$. 
 
 The {\em Girsanov condition} which makes the reweighting possible and whose validity we will check is
 \begin{equation}\label{Girs}
 \mathbb E \exp\left(M_s-\frac 12 \langle M\rangle_s\right) =1, \qquad s\ge 0,
 \end{equation}
 with $M_s := \int_0^s\left[\frac{\theta}{\sigma}-\frac{\sigma\gamma}{2}L_r(H_r)\right]dB_r$.
\begin{remark}\label{Girsa} If
 \eqref{Girs} holds, then there exists a measure $\tilde {\mathbb P}$ such that for all $s>0$
   $$\left. \frac{d\tilde{\mathbb P}}{d\mathbb P}\right|_{\mathcal{F}_s}   =
   \exp\left(M_s-\frac 12 \langle M\rangle_s\right),$$
   where $\mathcal F_s$ is the $\sigma$-field generated by $(H_r)_{0\le r\le s}$. By Girsanov's theorem (see Proposition \ref{Gir-BM} in the Appendix below), $\tilde B$ defined by \eqref{Btilde} then is a standard Brownian motion under $\tilde {\mathbb P}$. Thus, under $\tilde {\mathbb P}$, $H$ solves \eqref{reflecSDE} with $\tilde B$ instead of $B$. This gives the existence part of Proposition \ref{exun}, as soon as we have \eqref{Girs}.
 \end{remark}
 A sufficient condition for \eqref{Girs} is provided by the following lemma, which is Theorem~1.1, chapter 7 (page 152) in  \cite{Fr}.  \begin{lemma}\label{suffGirs}
  Assume that the quadratic variation of the continuous local martingale $M$ is of the form $\langle M\rangle_s = \int_0^s R_r \, dr$, and that for all $s>0$ there exist constants $a >0$ and $c<\infty$ such that   
  \begin{equation}\label{suffgir}
  \mathbb E \exp(aR_r) \le c,  \qquad 0\le r\le s.
  \end{equation}
  Then  \eqref{Girs} is satisfied.
 \end{lemma}
 In our situation, $R_r = \left| \frac{\theta}{\sigma}-\frac{\sigma\gamma}{2} L_r(H_r)\right|^2$, for which \eqref{suffgir} is implied by the following
 \begin{lemma}\label{Girstandard}
 Let $H$ be a Brownian motion   on $\mathbb R_+$ reflected at the origin, with variance parameter $v^2$.
Then for all $s>0$ there exists $\alpha = \alpha(s,v) >0$ and a constant $c < \infty$ such that   
$$\mathbb E\left(\exp (\alpha L_r(H_r)^2)\right)
\le c,   \qquad 0\le r\le s.$$
\end{lemma} 
\noindent
 {\em Proof:} Together with a simple scaling argument and a desintegration with respect to $H_r$, this is immediate from the following
\begin{lemma}\label{monloctime}
Let $\beta$ be a standard Brownian motion starting at $0$, and denote by $L_1(t)$ the local time accumulated by $|\beta|$ at position $t$ up to time $1$. There exist  constants $a >0$ and $c >0$ (not depending on $t$) such that 
 \begin{equation} \label{expmom}\mathbb E[e^{aL_1(t)^2}|\, |\beta_1| = t] \le c, \quad t\ge 0.
  \end{equation}
\end{lemma} \noindent
 {\em Proof:} Denote by $K_1(x)$ the local time of $\beta$ accumulated up to time $1$ at position $x$. First observe that for $t\ge 0$
 \begin{equation}\label{Lsum}
  L_1(t)= K_1(t) +  K_1(-t) \mbox{ a.s. }
 \end{equation}
For deriving \eqref{expmom}, by symmetry it suffices to condition under the event $\{\beta_1 = t\}$. Writing $\mathbb P^x$ for $\mathbb P[\, .\, | \beta_1= x]$ we conclude from \eqref{Lsum} and  the Cauchy-Schwarz inequality that for all $a>0$
\begin{equation}\label{CS}
\mathbb E^t[e^{aL_1(t)^2}] \le \left(\mathbb E^t\left[e^{4a K_1(t)^2}\right]\right)^{1/2} \left(\mathbb E^t\left[e^{4a K_1(-t)^2}\right]\right)^{1/2} .
\end{equation}
For $u \le 1$, the distribution of  $K_1(t)$ under $\mathbb P^t$ and conditioned under the event that $\beta$ hits $t$ first at time $u$, equals the  distribution of  $\sqrt{1-u}\, K_1(0)$ under $\mathbb P^0$. Similarly, for $u_1, u_2 \le 1$, the distribution of  $K_1(-t)$ under $\mathbb P^t$ and conditioned under the event that $(\beta_v)_{0\le v\le1}$ hits $-t$ first at time $u_1$ and last at time $u_2$, equals the the distribution of  $\sqrt{u_2-u_1}\, K_1(0)$ under $\mathbb P^0$. Consequently, 
\begin{equation}\label{compa}
\mathbb E^t\left[e^{4a K_1(t)^2}\right]\le \mathbb E^0\left[e^{4a K_1(0)^2}\right] , \qquad \mathbb E^t\left[e^{4a K_1(-t)^2}\right]\le \mathbb E^0\left[e^{4a K_1(0)^2}\right].
\end{equation}
By a result due to L\'evy (see formula (11) in \cite{Pi}),   $K_1(0)$ has under $\mathbb P^0$ a  Raleigh distribution, i.e.
$$\mathbb P^0(K_1(0)   > \ell ) = e^{-\frac 12 \ell^2}.$$
This means that $K_1^2(0)$ is  exponentially distributed, and hence, for suitably small $\delta>0$,  $\mathbb E^0\left[e^{\delta K_1(0)^2}\right]$ is finite. Now \eqref{expmom} follows from \eqref{CS} and \eqref{compa}. 
\epf
As stated above, Lemmas  \ref{suffGirs} and \ref{Girstandard} give together with Remark \ref{Girsa} the existence part of Proposition \ref{exun}.
For its uniqueness part, assume that $H$ is a weak solution of \eqref{reflecSDE}, governed by some measure $\mathbb P$.  For $n \in \mathbb N$, we define \begin{equation}\label{Tn}
T_n := \inf \{r > 0 : L_r(H_r) > n\}.
\end{equation} 
By a Girsanov transformation we can change the measure $\mathbb P$ into a measure $\bar {\mathbb P}$ under which, for all $n \in \mathbb N$,  
$$\bar B_{s\wedge n \wedge T_n} =  B_{s\wedge n \wedge T_n} + \int_0^{s\wedge n \wedge T_n} \left[\frac \theta \sigma -\gamma \frac {\sigma}2 L_r(H_r)\right]\, dr, \quad s\ge 0,$$
is a standard Brownian motion stopped at $n \wedge T_n$. Then, under  $\bar {\mathbb P}$, for all $n$, the process $H$ satisfies 
\begin{equation}\label{weaksolloc}
\left\{
\begin{aligned}
d  H_s &= \frac 2\sigma d \bar B_s+\frac 12 d  L_s(0), \quad 0 \le s \le n\wedge T_n, \\ \notag
H_0 &=0.
\end{aligned}
\right.
\end{equation}
Because the weak solution of \eqref{weaksolloc} is unique, the law of $(H_{s\wedge n \wedge T_n})_{s\ge 0}$ under $\bar {\mathbb P}$ is uniquely determined, and so is the law of $(H_{s\wedge n \wedge T_n})_{s\ge 0}$ under $\mathbb  P$, since it can be recovered by the inverse Girsanov transformation (note that $\mathbb P$ and  $\bar {\mathbb P}$ are mutually absolutely continuous on $\mathcal F_{n\wedge T_n}$).  From Theorem VI.1.7, page 225 of \cite{ReYo} one can infer that the local time $L=L_s(t)$ of $H$ is jointly continuous on $[0,\infty) \times [0,\infty)$ in its two arguments  $\mathbb P$-a.s. Consequently,  $T_n \to \infty$ as $n \to \infty$, $\mathbb P$-a.s, which shows weak uniqueness for \eqref{reflecSDE} and completes the proof of Proposition \ref{exun}.
\section{A discrete mass approximation}\label{discretemass}
The aim of this section is to set up a ``discrete mass - continuous time'' approximation of \eqref{fellog1} and \eqref{reflecSDE}. This will explain the intuition behind Theorem \ref{main}, and also will prepare for its proof. 

For $x> 0$ and $N\in \mathbb N$ the approximation of \eqref{fellog1} will be given by the  total mass $Z^{N,x}$  of a population of individuals, each of which has mass $1/N$. The initial mass is $Z_0^{N,x} =  \lfloor Nx \rfloor/N$, and $Z^{N,x}$ follows a Markovian jump dynamics: from its current state $k/N$, 
\begin{eqnarray}\label{Zndyn}
Z^{N,x} \mbox{ jumps to } \begin{cases}
(k+1)/N \mbox{ at rate } kN\sigma^2/2 + k\theta\\ (k-1)/N \mbox{ at rate }  kN\sigma^2/2 + k(k-1)\gamma/N. \end{cases} 
\end{eqnarray}
For $\gamma =0$, this is (up to the mass factor $1/N$) as a Galton-Watson process in continuous time: each individual independently spawns a child at rate $N\sigma^2/2 + \theta$, and dies (childless) at rate $N\sigma^2/2$. For $\gamma \neq 0$, the additional quadratic death rate destroys the independence, and hence also the branching property. However, when viewing the individuals alive at time $t$ as being arranged ``from left to right'', and by decreeing that each of the pairwise fights (which happen at rate $2\gamma$ and always end lethal for one of the two involved individuals) is won by the individual to the left, we arrive at the additional death rate $2\gamma {\mathcal L}_i(t)/N$ for individual $i$, where ${\mathcal L}_i(t)$ denotes the number of individuals  living at time $t$ to the left of individual $i$.

 The just described reproduction dynamics gives rise to a {\em forest} $F^{N,x}$ of {\em trees of descent}, drawn into the plane as sketched in Figure 1. At any branch point, we imagine the ``new'' branch being placed to the right of the mother branch. Because of the asymmetric killing, the trees further to the right have a tendency to stay smaller: they are ``under attack'' by the trees to their left. Note also that, with the above described construction, the $F^{N,x}$, $x>0$, are coupled: when $x$ is increased by $1/N$, a new tree is added to the right. We denote the union of the $F^{N,x}$, $x>0$, by $F^N$.

From $F^N$ we read off a continuous and piecewise linear $\mathbb R_+$-valued path $H^N= (H^N_s)$ (called the {\em exploration path} of $F^N$) in the following way:

Starting from the root of the leftmost tree, one goes upwards at speed $2N$ until one hits the top of the first mother branch (this is the leaf marked with $\mathfrak D$ in Figure 1). There one turns and goes downwards, again at speed $2N$, until arriving at the next branch point (which is  $\mathfrak B$ in Figure 1). From there one goes upwards into the (yet unexplored) next branch, and proceeds in a similar fashion until  being back at height $0$, which means that the exploration of the leftmost tree is completed. Then explore the next tree, and so on.
\begin{figure}\label{figure1}
\psfrag{t}{$t$}
\psfrag{s}{$s$}
\psfrag{F}{$F^N$}
\psfrag{H}{$H^N$}
\psfrag{L}{$\mathfrak D$}
\psfrag{B}{$\mathfrak B$}
\includegraphics[width=15cm]{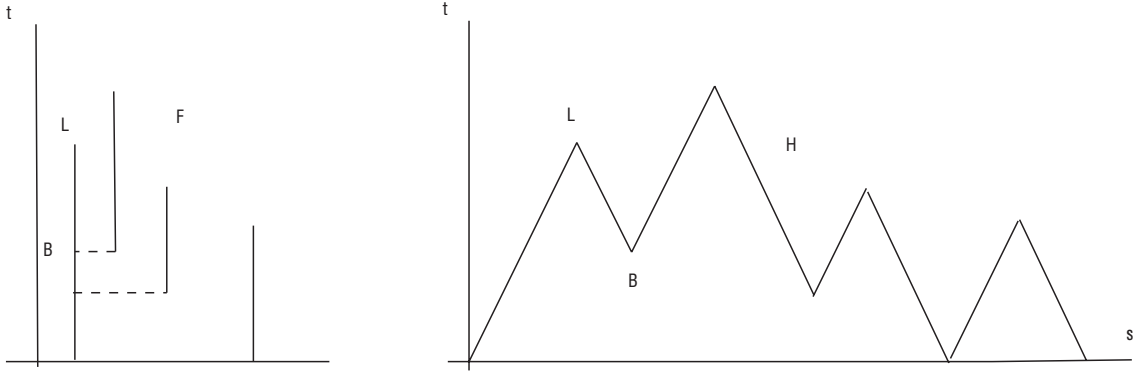}
\caption{A realization of (the first two trees of) $F^N$ and (the first two excursions of) its exploration 
path $H^N$. The $t$-axis is {\em real time} as well as {\em exploration height}, the 
$s$-axis is {\em exploration time}.}
\end{figure}
For $x >0$ we denote by $S^N_x$ the time at which the exploration of the forest $F^{N,x}$ is completed.
Obviously, for each $t \ge 0$, the number of  branches in $F^{N,x}$ that are alive at time $t$ equals half the number of $t$-crossings of the  exploration path of $F^N$ stopped at $S^N_x$. Recalling that the slope of $H^N$ is $\pm 2N$, we define
\begin{equation}\label{defLN}
\Lambda_s^N(t) := \frac 1{2N} \# \mbox{ of $t$-crossings of $H^N$ between exploration times $0$ and $s$,}
\end{equation}
where we count a local minimum of $H^N$ at $t$ as two $t$-crossings, and a local maximum as none. Note that by our convention both $s\mapsto \Lambda_s^N(t)$ and $t\mapsto  \Lambda_s^N(t)$ are right continuous, and in particular  $\Lambda_0^N(0) = 0$.  We call $\Lambda_s^N(t)$ the  {\em (unscaled) local time} of $H^N$ accumulated at height $t$ up to time $s$. This name is justified also by the following {\em occupation times formula}, valid for all measurable $f:\mathbb R_+ \to \mathbb R_+$
\begin{equation}\label{occtimeN}
\int_0^s f(H_r^N)\, dr = \int_0^\infty f(t)\, \Lambda_s^N(t)\, dt, \qquad s\ge 0.
\end{equation} 
The exploration time $S^N_x$ which it takes to traverse all of the $ \lfloor Nx\rfloor$ trees in $F^{N,x}$ can be expressed as \begin{equation}\label{SN}
S^N_x = \inf \{s : \Lambda_s^N (0) \ge \lfloor Nx\rfloor/N\}. 
\end{equation}
\begin{proposition} \label{dynHN}
The exploration path $s\mapsto H^N_s$ obeys the following stochastic dynamics:
\begin{itemize}
\item At time $s=0$, $H^N$ starts at height $0$ and with slope $2N$. 
\item
While $H^N$ moves upwards, its slope jumps from $2N$ to $-2N$ at rate $N^2\sigma^2+ 4\gamma N \ell$, where $\ell=\Lambda_s^N(H^N_s)$ is the local time accumulated by $H^N$ at the current height $H_s^N$ up to the current exploration time $s$.
\item
While $H^N$ moves downwards, its slope jumps from $-2N$ to $2N$ at rate $N^2\sigma^2+ 2N\theta$.
\item
Whenever $H^N$ reaches height $0$, it is reflected above $0$.
\end{itemize}
\end{proposition} \noindent
\bpf We give here an informal proof  which contains the essential ideas. (A more formal proof can be carried out along the arguments of the proof of Theorem 2.4 in \cite{BPS}.)

Recall that the death rate of an individual $i$ living at real time $t$ is  $ N\sigma^2/2+ 2\gamma {\mathcal L}_i(t)/N$, where   ${\mathcal L}_i(t)$ is the number of individuals living at time $t$ to the left of individual $i$. Assume the individual $i$ living at time $t$ is explored first at time $s$, hence $H_s = t$, and $H$ has slope $2N$ at time $s$. Because of \eqref{defLN}, while $H^N_s$ goes upward, we have ${\mathcal L}_i(t) = N\Lambda_s(H_s^N)$. The rate in $t$ is the rate in $s$  multiplied by the factor $2N$ which is the absolute value of the slope. This gives the claimed jump rate $2N(N\sigma^2/2+ 2\gamma \Lambda_s(H_s^N))$ from slope $2N$ to slope $-2N$, which can be seen as the rate at which  the ``death clock'' rings (and leads to a downward jump of the slope) along the rising pieces of the exploration path $H^N$. On the other hand, the ``birth clock'' rings along the falling pieces of $H^N$, its rate being 
$ N\sigma^2/2+\theta$ in real time and $2N(N\sigma^2/2+\theta)$ in exploration time, as claimed in the proposition. Note that the process of birth times along an individual's lifetime is a homogeneous Poisson process which (in distribution) can as well be run backwards from the individual's death time. Also note that, due to the ``depth-first-search''-construction of $H^N$, along falling pieces of $H^N$ always yet unexplored parts of the forest are visited as far as the birth points are concerned.
\epf

The next statement is a discrete version of Theorem \ref{main}, and will later be used for the proof of Theorem \ref{main} by taking $N \to \infty$.  
\begin{corollary}\label{discreteRN}
Let $H^N$ be a stochastic process following the dynamics specified in Proposition \ref{dynHN}, and $\Lambda^N$ be its local time as defined by \eqref{defLN}. For $x > 0$, let $S_x^N$ be the stopping time defined by \eqref{SN}. Then $t\mapsto \Lambda^N_{S^N_x}(t)$ follows the jump dynamics \eqref{Zndyn}.
\end{corollary}
\bpf
By Proposition \ref{dynHN}, $H^N$ is equal in distribution to the exploration path of the random forest $F^N$. Hence $\Lambda_{S^N_x}(t)$ is equal in distribution to $Z^{N,x}_t$, where $NZ^{N,x}_t$ is the number of branches alive in $F^{N,x}$ at time $t$. Since $Z^{N,x}$ follows the dynamics  \eqref{Zndyn}, so does $\Lambda^N_{S^N_x}$.
\epf
The next lemma will also be important in the proof of Theorem \ref{main}.
\begin{lemma}\label{extendable}
Let $H^N$ and $S^N_x$ be as in Corollary \ref{discreteRN}. Then $S^N_x \to \infty$ a.s. as $x \to \infty$.
\end{lemma}
\bpf
Consider $x = a/N$ for $a=1,2\ldots$. Applying \eqref{occtimeN} with $s = S^N_x$ and $f\equiv 1$ we obtain the equality $S_x^N = \int_0^\infty \Lambda^N_{S^N_x}(t)\, dt$.

 According to Corollary \ref{discreteRN}, $\Lambda^N_{S^N_x}(t)$ follows the jump dynamics $\eqref{Zndyn}$, with initial condition $\Lambda_{S^N_x}(t)=a/N$. By coupling $\Lambda^N_{S^N_x}(t)$ with a ``pure death process'' $K^N$ that starts in $a/N$ and jumps from $k/N$ to $(k-1)/N$ at rate $k(k-1)(N\sigma^2/2 + \gamma/N)$, we see that $\int_0^\infty \Lambda^N_{S^N_x}(t)\, dt$ is stochastically bounded from below by $\int_0^{T_2}K^N_tdt$, where $T_2$ is the first time at which $K^N$ comes down to $2/N$. The latter integral equals a sum of independent exponentially distributed random variables with parameters $(j-1)(N^2\sigma^2/2+\gamma))$, $j=2,\ldots, a$. This sum diverges as $a\to \infty$.
 \epf

\section{Convergence of the mass processes $Z^{N,x}$ as $N\to \infty$}\label{secZ}
The process $\{Z^{N,x}_t,\ t\ge0\}$ with dynamics \eqref{Zndyn} is a Markov process with values in the set
$E_N:=\{k/N,\ k\ge1\}$, starting from $\lfloor Nx\rfloor/N$, with generator $A^N$ given by 
\begin{align}
& A^Nf(z)=Nz \left( N\frac{\sigma^2}{2}+\theta \right) \left[f\left(z+\frac{1}{N}\right)-f(z)\right]
\\ \notag  & \phantom{aaaaaaaaaaaaaaaaaa} +Nz\left(N\frac{\sigma^2}{2}
+\gamma \left(z-\frac{1}{N}\right)\right) \left[f\left(z-\frac{1}{N}\right)-f(z)\right],
\end{align}
for any $f:E_N\to\R$, $z\in E_N$. (Note that the distinction between symmetric and ordered killing is irrelevant here.) Applying successively the above formula to the cases 
$f(z)=z$ and $f(z)=z^2$, we get that
\begin{align}
Z^{N,x}_t&=Z^{N,x}_0+\int_0^t\left[\theta Z^{N,x}_r-\gamma Z^{N,x}_r\left(Z^{N,x}_r-\frac{1}{N}\right)\right]dr+M^{(1)}_t,\label{z}\\
\left(Z^{N,x}_t\right)^2&= \left(Z^{N,x}_0\right)^2 +
2\int_0^tZ^{N,x}_r\left[\theta Z^{N,x}_r-\gamma Z^{N,x}_r\left(Z^{N,x}_r-\frac{1}{N}\right)\right]dr
\nonumber\\
&\qquad +\int_0^t\left[\sigma^2Z^{N,x}_r+\frac{\theta }{N}Z^{N,x}_r+\frac{\gamma}{N}\left(Z^{N,x}_r-\frac{1}{N}\right)Z^{N,x}_r\right]dr+M^{(2)}_t,\label{zsquare}
\end{align}
where $\{M^{(1)}_t,\ t\ge0\}$ and $\{M^{(2)}_t,\ t\ge0\}$ are local martingales. It follows from \eqref{z} and \eqref{zsquare} that
\begin{equation}\label{bracket}
\langle M^{(1)}\rangle_t=
\int_0^t\left[\sigma^2Z^{N,x}_r+\frac{\theta }{N}Z^{N,x}_r+\frac{\gamma}{N}\left(Z^{N,x}_r-\frac{1}{N}\right)Z^{N,x}_r\right]dr.
\end{equation}
We now prove
\begin{lemma}\label{4thmoment} 
For any $T>0$,
$$\sup_{N\ge1}\sup_{0\le t\le T}\E\left[\left(Z^{N,x}_t\right)^4\right]<\infty.$$
\end{lemma}
\noindent
An immediate Corollary of this Lemma is that $\{M^{(1)}_t\}$ and $\{M^{(2)}_t\}$ are 
in fact martingales.

\bpf The same computation as above, but now with $f(z)=z^4$, gives
\begin{equation}\label{4moment}
\left(Z^{N,x}_t\right)^4= \left(Z^{N,x}_0\right)^4 +\int_0^t\Phi_N\left(Z^{N,x}_r\right)dr
+M^{(4)}_t,
\end{equation}
where $\{M^{(4)}_t,\ t\ge0\}$ is a local martingale
and for some $c>0$ independent of $N$,
\begin{equation}\label{ineq4}
\Phi_N(z)\le c(1+z^4).
\end{equation}
We note that $NZ^{N,x}_t$ is bounded by the value at time $t$ of a Yule process
(which arises when suppressing the deaths), which is a finite sum of mutually independent  geometric random variables, hence has
finite moments of any order. Hence $M^{(4)}$ is in fact a martingale.
We then can take the expectation in \eqref{4moment}, and deduce from \eqref{ineq4} and Gronwall's Lemma that
for $0\le t\le T$,
$$\E\left[\left(Z^{N,x}_{t}\right)^4\right]\le \left[\left(Z^{N,x}_0\right)^4+cT\right] e^{cT},$$ which  implies the result.
\epf
We shall also need below the
\begin{lemma}\label{supZ}
For any $T>0$,
$$\sup_{N\ge1}\E\left[\sup_{0\le t\le T}\left(Z^{N,x}_t\right)^2\right]<\infty.$$
\end{lemma}
\bpf Since from \eqref{z}, $Z^{N,x}_t\le Z^{N,x}_0+\theta\int_0^t Z^{N,x}_rdr+M^{(1)}_t$,
$$
\sup_{r\le t}|Z^{N,x}_r|^2\le3 |Z^{N,x}_0|^2+3 t\theta^2\int_0^t|Z^{N,x}_r|^2dr
+3\sup_{r\le t}|M^{(1)}_r|^2.$$
This together with \eqref{bracket}, Doob's $L^2$-inequality for martingales and Lemma \ref{4thmoment}
implies the result. 
\epf

Remark \ref{tightnessCriteria} in the Appendix combined with \eqref{z}, \eqref{bracket} and Lemma \ref{4thmoment} guarantees that the tightness of $\{Z^N_0\}_{n\ge1}$ implies that of $\{Z^N\}_{N\ge1}$ in $D([0,+\infty))$.

Standard arguments exploiting \eqref{z} and \eqref{zsquare} now allow us to deduce the convergence of the mass processes (for a detailed proof, see e.g. Theorem 5.3 p. 23 in \cite{SM}).
\begin{proposition}\label{convZ}
As $N\to\infty$, $Z^{N,x}\Rightarrow Z^x$,
where $Z^x$ is the unique solution of the SDE \eqref{fellog1}
and thus is a Feller diffusion with logistic growth.
\end{proposition} 

\section{Convergence of the exploration path in the case $\theta = \gamma = 0$.}\label{allzero}
Let $H^N$ be a stochastic process as in Proposition \ref{dynHN} with $\theta = \gamma = 0$. The aim of this section is to provide a version of the joint convergence (as $N\to \infty$) of $H^N$ and its local time  which is suitable for the change of measure that will be carried through in Section \ref{secgirs}. This is achieved in Theorem \ref{conHL} and its Corollary \ref{convH-LS}. The proof of Theorem 2 is carried out in two major parts. The first part (Proposition \ref{convcombin}) provides a refined version of the joint convergence of $H^N$ and its local time at level $0$, the second part (starting from Lemma \ref{contLT}) extends this to the other levels as well.

We define the {\em (scaled) local time accumulated by $H^N$ at level $t$ up to time $s$} as
\begin{align}\notag
  L^N_s(t):= \frac{4}{\sigma^2}\lim_{\eps\to0}\frac{1}{\eps}\int_0^{s}{\bf1}_{\{t< H^N_u<t+\eps\}}du
\end{align}
Note that this process is neither right-- nor left--continuous as a function of $s$. However since the jumps are of size
  $O(1/N)$, the limit of $L^N$ as $N\to\infty$ will turn out to be continuous. In fact,
we will show that $L^N$ converges as $N\to \infty$ to the semi--martingale local time of the limiting process $H$, hence the scaling factor $4/\sigma^2$. It is readily checked that  $L_{0+}^N(0) = \frac 2{N\sigma^2}$, and 
\begin{equation}\label{corresp}
L_s^N(t)= \frac 4{\sigma ^2} \Lambda_s^N(t), \quad \forall s\geq 0, t\geq 0,
\end{equation}
where $\Lambda$ was defined in \eqref{defLN}. Then with $S_x^N$ defined in \eqref{SN}, we may rewrite
\begin{equation}\label{defSN}
S_x^N=\inf \{s>0: L_s^N(0)\geq  \frac4{\sigma ^2}\lfloor Nx\rfloor/N\}.
\end{equation}
From \eqref{corresp} and Corollary \ref{discreteRN} we see that
\begin{equation*}
Z_t^{N,x}:= \frac{\sigma ^2}{4} L_{S_x^N}^N(t)
\end{equation*}
follows the jump dynamics \eqref{Zndyn} in the case $\theta =\gamma =0$. 

Let $\{V_s^N, s\geq 0\}$ be the c\`adl\`ag $\{-1,1\}$-valued process which is such that for a.\ a. $s>0$,
\begin{equation*}
\frac{dH_s^N}{ds}= 2N V_s^N.
\end{equation*}
We can express $L^N$ in terms of $H^N$ and $V^N$ as
\begin{align}\notag
  L^N_s(t)= \frac{4}{\sigma^2}\frac{1}{2N}\left[\sum_{0\le r<s}{\bf1}_{\{H^N_r=t\}}
  \left(1+\frac{1}{2}(V_r^N-V_{r-}^N)\right)+{\bf1}_{\{H^N_s=t\}}{\bf1}_{\{V^N_{s-}=-1\}}\right]
  \end{align}
  where we put $V_{0-}=+1$. Note that 
  any $r<s$ at which $t$ is a local minimum of $H^N$ counts twice in the sum of the last line, while any
  $r<s$ at which $t$ is a local maximum of $H^N$ is not counted in the sum.

We have
\begin{equation}\label{HnVn}
\begin{split}
H_s^N&= 2N\int_0^sV_r^Ndr \, , \\
V_s^N&= 1+ 2\int _0^s \mathbf{1}_{\{V_{r-}^N= -1\}} dP^N_r - 2\int _0^s \mathbf{1}_{\{V_{r-}^N= 1\}} dP^N_r + \frac {N\sigma ^2}2(L_{s}^N(0)-L_{0+}^N(0)),
\end{split}
\end{equation}
where $ \{P^N_s, s\geq 0\}$ is a Poisson process with intensity $N^2\sigma ^2$.  Note that $M_s^N= P^N_s- N^2\sigma ^2s,$ $ s\geq 0$, is a martingale. 
The second equation is prescribed by the statement of Proposition \ref{dynHN} (in fact its simplified version in case $\theta=\gamma=0$)~: the initial velocity is positive, whenever $V^N=1$, it jumps to -1 (i. e. makes a jump of size -2) at rate $N^2\sigma ^2$, whenever $V^N=-1$, it jumps to +1 (i. e. makes a jump of size +2) at rate $N^2\sigma ^2$, and in addition it makes a jump of size 2 whenever $H^N$ hits 0.
Writing $V^N_r$ in the first line of \eqref{HnVn} as 
$$\mathbf{1}_{\{V_{r-}^N= +1\}}-\mathbf{1}_{\{V_{r-}^N= -1\}}$$
and denoting by $M^N_s$ the martingale $P^N_s-N^2\sigma^2s$, we deduce from 
\eqref{HnVn}
\begin{equation}\label{HN}
H_s^N+ \frac{V_s^N}{N\sigma ^2}=  M_s^{1,N}- M_s^{2,N}+\frac 12 L_{s}^N(0) -\frac 12 L_{0^+}^N(0),
\end{equation}
where
\begin{equation}\label{twomart}
M_s^{1,N}= \frac2{N\sigma ^2} \int _0^s \mathbf{1}_{\{V_{r-}^N= -1\}} dM_r^N \quad \mbox{ and } \quad M_s^{2,N}= \frac2{N\sigma ^2} \int _0^s \mathbf{1}_{\{V_{r-}^N= 1\}} dM_r^N
\end{equation}
are two mutually orthogonal martingales. Thanks to an averaging property of the $V^N$ (see step 2 in the proof of Proposition  \ref{convcombin} below) these two martingales will converge as $N\to \infty$ to two independent Brownian motions with variance parameter $2/\sigma^2$ each. Together with the appropriate convergence of $L^N(0)$, \eqref{HN} then gives the required convergence of $H^N$, see Proposition \ref{convcombin}. We are now ready to state the main result of this section.
\begin{theorem}\label{conHL} For any $x>0$, as $N \to \infty$,
\begin{equation*}
(\{H_s^N, M^{1,N}_s,M^{2,N}_s,\ s\geq 0\},\{L_s^N(t), s,t\geq 0\}, S^N_x)\Rightarrow (\{H_s, \frac{\sqrt{2}}{\sigma}B^1_s,
\frac{\sqrt{2}}{\sigma}B^2_s,\ s\geq 0\},\{L_s(t), s,t\geq 0\},S_x)
\end{equation*}
for the topology of locally uniform convergence in $s$ and $t$. $B^1$ and $B^2$ are two mutually independent standard Brownian motions,  $H$ solves the SDE \eqref{BR} whose driving Brownian motion $B$ is given as
$$B_s=\frac{1}{\sqrt{2}}(B^1_s-B^2_s),$$  $L$ is the semi--martingale local time of $H$, and $S_x$ has been defined in \eqref{Sx}.
\end{theorem}

An immediate consequence of this result is
\begin{corollary}\label{convH-LS}
For any $x>0$, as $N\to\infty$,
\begin{equation*}
(\{H_s^N, M^{1,N}_s,M^{2,N}_s,\  s\geq 0\},\{L_{S^N_x}^N(t), t\geq 0\})\Rightarrow (\{H_s, \frac{\sqrt{2}}{\sigma}B^1_s,
\frac{\sqrt{2}}{\sigma}B^2_s,\  s\geq 0\},\{L_{S_x}(t), t\geq 0\})
\end{equation*}
in $C([0,\infty))\times (D([0,\infty)))^3$.
\end{corollary}
Recall (see Lemma \ref{CdansD} in the Appendix) that convergence in $D([0,\infty))$ is equivalent to locally uniform convergence, provided the limit is continuous. Also note that in the absence of reflection, the weak convergence of $H^N$ to Brownian motion would be a consequence of Theorem 7.1.4 in \cite {EK}, and would be a variant of ``Rayleigh's random flight model'', see Corollary 3.3.25 in \cite{DWS}.

A first preparation for the proof of Theorem \ref{conHL} is 
\begin{lemma}\label{HNtight} The sequence $\{H^N\}$ is tight in $C([0,\infty))$.
\end{lemma}
\bpf
To get rid of the local time term in \eqref{HN}, we consider a process $R^N$ of which $H^N$ is the absolute value. More explicitly, let $(R^N,W^N)$ be the $\R\times\{-1,1\}$--valued process that solves the system (which is exactly \eqref{HnVn} without reflection)
\begin{align*}
R_s^N&= 2N \int_0^s W_r^N dr \, ,\\
W_s^N&= 1+ 2\int _0^s  \mathbf{1}_{\{W_{r-}^N= -1\}} dP^N_r - 2\int _0^s  \mathbf{1}_{\{W_{r-}^N= +1\}} dP^N_r.
\end{align*}
We observe that
\begin{equation*}
(H^N,V^N) \equiv (|R^N|,\,  \mathrm{sgn}(R^N) W^N)\, .
\end{equation*}
Clearly tightness of $\{R^N\}$ will imply that of $\{H^N\}$, since $|H_s^N- H_t^N| \leq |R_s^N- R_t^N|$ for all $s, t \ge 0$.
Now we have 
\begin{equation*}
R_s^N+ \frac{W_s^N}{N\sigma ^2}= \frac 1{N \sigma ^2} - \frac 2{N \sigma ^2} \int_0^s W_{r-}^N \, dM_r^N.
\end{equation*}
By Proposition \ref{sc} in the Appendix, the sequence $\{R^N\}_{N\geq 1}$ is tight, and so is $\{H^N\}_{N\geq 1}$.
\epf
\begin{proposition}\label{convcombin} 
Fix $x>0$.
As $N\to\infty$,
\begin{align*}
\big( H^N, M^{1,N}, M^{2,N}, L^N(0), S_x^N \big) \Rightarrow  &\Big( H, \frac{\sqrt{2}}{\sigma }B^1, \frac{\sqrt{2}}{\sigma }B^2, L(0), S_x\Big) \\
& in \quad C([0,\infty )) \times \big(D([0,\infty ))\big)^3 \times [0,\infty ),
\end{align*}
where $B^1$ and $B^2$ are two mutually independent standard Brownian motions, and  $H$ solves the  SDE
\begin{equation}\label{SDH}
H_s= \frac2{\sigma }B_s + \frac 12 L_s(0), \quad s\ge 0,
\end{equation}
with $B_s := (1/\sqrt 2) (B_s^1- B_s^2)$, and $L(0)$ denoting the local time at level $0$ of $H$. (Note that $B$ is again a standard Brownian motion.)
\end{proposition} \bpf 
The proof is organized as follows. Step 1 establishes the weak convergence of 
$(H^N,M^{1,N},M^{2,N},L^N(0))$ along a subsequence. Step 2 and step 3 together characterize the law of the limiting two--dimensional martingale, step 4 identifies the limit of the local time term. In step 5 we note that the entire sequence converges. Finally step 6 takes the limit in the quintuple (including $S_x^N$).

\bigskip\noindent
{\sc Step 1.}
Note that 
\begin{description}
\item[\rm{i)}] from Lemma \ref{HNtight}, the sequence $\{H^N, N\geq 1\}$ is tight in $C([0,\infty ])$;
\item[\rm{ii)}]$\sup_{s\geq 0} \frac{|V_s^N|}{N\sigma ^2}\rightarrow 0$ in probability as $N\rightarrow \infty $;
\item[\rm{iii)}]from Proposition \ref{sc}, $\{M^{1,N}, N\geq 1\}$ and $\{M^{2,N}, N\geq 1\}$ are tight in $D([0,\infty ])$, any limiting martingales $M^1$ and $M^2$ being continuous;
\item[\rm{iv)}]it  follows from the first $3$ items, \eqref{HN} and Proposition \ref{tightcombined} that $\{L_s^N(0), N\geq 1\}$ is tight in $D([0,\infty ])$, the limit $K$ of any converging subsequence being continuous and increasing.
\end{description}
Working along a diagonal subsequence we can extract a subsequence, still denoted as an abuse like the original sequence, such that along that subsequence
\begin{equation*}
\big( H^N, M^{1,N}, M^{2,N}, L^N(0) \big) \Rightarrow  \big( H, M^1, M^2, K \big).
\end{equation*}
\bigskip
\noindent
{\sc Step 2.} We claim that for any $s> 0$,
\begin{equation*}
 \int_0^s \mathbf{1}_{\{V_{r}^N= 1\}} dr \rightarrow \frac{s}2,\quad \int_0^s \mathbf{1}_{\{V_{r}^N= -1\}} dr \rightarrow \frac{s}2
\end{equation*}
in probability, as $N\rightarrow \infty $. 
This follows by taking the limit in the sum and the difference of the two following identities~:
\begin{align*}
\int_0^s \mathbf{1}_{\{V_{r}^N= 1\}} dr+\int_0^s \mathbf{1}_{\{V_{r}^N= -1\}} dr&=s,\\
\int_0^s \mathbf{1}_{\{V_{r}^N= 1\}} dr-\int_0^s \mathbf{1}_{\{V_{r}^N= -1\}} dr&=(2N)^{-1}H^N_s,
\end{align*}
since $H^N_s/N\to0 $ in probability, as $N\to\infty$, thanks to Lemma \ref{HNtight}.

\bigskip
\noindent{\sc Step 3.} By Step 1 iii), $M_s=\begin{pmatrix}M^1_s\\ M^2_s\end{pmatrix}$, the weak limit of $\begin{pmatrix}M^{1,N}_s\\ M^{2,N}_s\end{pmatrix}$ along the chosen subsequence, is a 2--dimensional continuous martingale. In order to identify it, we first introduce some useful notation.  We write $M_s^{\otimes2}$ for the $2\times 2$ matrix whose $(i,j)$--entry equals $M_s^i\times M_s^j$, and
 $\langle\langle M\rangle\rangle_s$ for the $2\times 2$ matrix--valued predictable increasing process which is such that
$$M_s^{\otimes2}-\langle\langle M\rangle\rangle_s$$ is a martingale, and note that the $(i,j)$--entry of the matrix $\langle\langle M\rangle\rangle_s$ equals $\langle M^i,M^j\rangle_s$. We adopt similar notations for the pair $M^{1,N}$, $M^{2,N}$.

From Step 2 we deduce that, as $N\to\infty$,
\begin{align*}
\langle\langle \begin{pmatrix}M^{1,N}\\ M^{2,N}\end{pmatrix} \rangle\rangle_s 
&=  \frac{4}{\sigma ^2} \int_0^s \begin{pmatrix}\mathbf{1}_{\{V_{r}^N= -1\}}&0\\ 
0& \mathbf{1}_{\{V_{r}^N= 1\}}\end{pmatrix}dr \\
&\to \frac{2}{\sigma ^2}s I
\end{align*}
in probability, locally uniformly in $s$, where $I$ denotes the $2\times2$ identity matrix. Consequently
$$\begin{pmatrix}M^{1,N}_s\\ M^{2,N}_s\end{pmatrix}^{\otimes2}- 
\langle\langle \begin{pmatrix}M^{1,N}\\ M^{2,N}\end{pmatrix} \rangle\rangle_s\Rightarrow
\begin{pmatrix}M^{1}_s\\ M^{2}_s\end{pmatrix}^{\otimes2}-\frac2{\sigma ^2}sI$$
in $D([0,\infty);\R^4)$ as $N\to\infty$, and since the weak limit of martingales is a local martingale, there exist two mutually independent standard Brownian motions $B^1$ and $B^2$ such that
\begin{equation*}
M_s^1= \frac{\sqrt{2}}{\sigma } B_s^1, \, M_s^2= \frac{\sqrt{2}}{\sigma } B_s^2,\  s\geq 0.
\end{equation*}
Taking the weak limit in \eqref{HN} we deduce that
\begin{align*}
H_s&=\frac{\sqrt{2}}{\sigma}(B^1_s-B^2_s)+\frac{1}{2}K_s\\
&=\frac{\sqrt{2}}{\sigma}B_s+\frac{1}{2}K_s,
\end{align*}
where $B_s=(B^1_s-B^2_s)/\sqrt{2}$ is also a standard Brownian motion.

\bigskip\noindent
{\sc Step 4.} 
For each $\ell \geq 1$, we define the function $f_{\ell}: \mathbb{R}_{+} \rightarrow [0,1]$ by $f_{\ell} (x)= (1-\ell x)^{+}$. We have that for each $N,\ell\geq 1, s>0$, since $L^N(0)$ increases only when $H^N= 0$,
\begin{equation*}
\mathbb{E} \Big( \int_0^s f_\ell (H_r^N) dL_r^N(0)- L_s^N(0) \Big) \geq 0.
\end{equation*}
Thanks to Lemma \ref{convStieltjes} in the Appendix we can take the limit in this last inequality as $N \rightarrow \infty $, yielding
\begin{equation*}
\mathbb{E} \Big( \int_0^s f_\ell (H_r) dK_r - K_s \Big) \geq 0.
\end{equation*}
Then taking the limit as $\ell \rightarrow \infty $ yields
\begin{equation*}
\mathbb{E} \Big( \int_0^s  \mathbf{1}_{\{H_r= 0\}} dK_r - K_s \Big) \geq 0.
\end{equation*}
But the random variable under the expectation is clearly nonpositive, hence it is zero a.s., in other words
\begin{equation*}
K_s= \int_0^s  \mathbf{1}_{\{H_r= 0\}} dK_r, \quad \forall s\geq 0,
\end{equation*}
which means that the process $K$ increases only when $H_r=0$.

From the occupation times formula
\begin{equation*}
\frac4{\sigma ^2} \int_0^s g(H_r) dr = \int_0^{\infty } g(t) L_s(t) dt
\end{equation*}
applied to the function $g(h)= \mathbf{1}_{\{h=0\}}$, we deduce that the time spent by the process $H$ at $0$ has a.s. zero Lebesgue measure. Consequently
\begin{equation*}
\int_0^s  \mathbf{1}_{\{H_r= 0\}} dB_r \equiv 0 \quad a.s.
\end{equation*}
hence a.s.
\begin{equation*}
B_s= \int_0^s  \mathbf{1}_{\{H_r> 0\}} dB_r \quad \forall s\geq 0.
\end{equation*}
It then follows from Tanaka's formula applied to the process $H$ and the function $h\rightarrow h^{+}$ that $K=L(0)$. 

\bigskip
\noindent {\sc Step 5.} We have proved so far  that $Q^N\Rightarrow Q$ along some subsequence, where
$Q^N=(H^N, M^{1,N}, M^{2,N}, L^N(0))$, $Q=(H, \frac{\sqrt{2}}{\sigma }B^1, \frac{\sqrt{2}}{\sigma }B^2, L(0))$. 
Note that not only subsequences but the entire sequence $Q^1, Q^2, Q^3, \ldots$ converges, since the limit law 
is uniquely characterized.

\bigskip
\noindent
{\sc Step 6.} It remains to check that for any $x>0$, as $N\to\infty$,
\begin{equation*}
(Q^N, S_x^N) \Rightarrow (Q, S_x) \quad 
\text {in } C([0,\infty ]) \times \big(D([0,\infty ])\big)^3 \times [0,\infty ].
\end{equation*}
To this end, let us define the function $\Phi $ from $\mathbb{R}_{+}\times C_{\uparrow }(\mathbb{R}_{+})$ into $\mathbb{R}_{+}$ by
\begin{equation*}
\Phi (x,y)= \inf \{s>0: y(s)> \frac4{\sigma ^2}x\}.
\end{equation*}
For any $x$ fixed, the function $\Phi (x,.)$ is continuous in the neighborhood of a function $y$ which is strictly increasing at the time when it first reaches the value $x$. Clearly $S_x= \Phi (x, L_{.}(0))$. Define
\begin{equation*}
S_x^{'N}:=  \Phi (x, L_{.}^N(0)).
\end{equation*}
We note that for any $x>0, s\mapsto L_s(0)$ is a.s. strictly increasing at time $S_x$, which is a stopping time. This fact follows from the strong Markov property, the fact that $H_{S_x}=0$, and $L_{\varepsilon }(0)>0$, for all $\varepsilon >0$. Consequently $S_x$ is a.s. a continuous function of the trajectory $L_{.}(0)$ , then also of $Q$, and
\begin{equation*}
(Q^N, S_x^{'N}) \Rightarrow (Q, S_x).
\end{equation*}
It remains to prove that $S_x^{'N}- S_x^N\rightarrow 0$ in probability. For any $y<x$ and $N$ large enough
\begin{equation*}
0\leq S_x^{'N}- S_x^N \leq S_x^{'N}- S_y^{'N}.
\end{equation*}
Clearly $S_x^{'N}- S_y^{'N}\Rightarrow S_x-S_y$, hence for any $\varepsilon >0$,
\begin{equation*}
0\leq \limsup_{N} \mathbb{P}(S_x^{'N}- S_x^N\geq \varepsilon )\leq  \limsup_{N} \mathbb{P}(S_x^{'N}- S_y^{'N}\geq \varepsilon )\leq  \mathbb{P}(S_x- S_y\geq \varepsilon ).
\end{equation*}
The result follows, since $S_y\rightarrow S_{x-}$ as $y\rightarrow x,y<x$, and $S_{x-}=S_x$ a.s.
\epf

For the proof of Theorem \ref{conHL} we will need the following lemmata:
\begin{lemma}\label{naive}
For any $s>0$, $t>0$, the following identities hold a.s.
\begin{align*}
(H^N_s-t)^+&=2N\int_0^sV^N_r{\bf1}_{\{H^N_r>t\}}dr,\\
V^N_s{\bf1}_{\{H^N_s>t\}}&=\frac{\sigma^2 N}{2}L^N_s(t)
+\int_0^s{\bf1}_{\{H^N_r>t\}}dV^N_r.
\end{align*}
\end{lemma}
\bpf
The first identity is elementary, and is true along any piecewise linear, continuous trajectory $\{H^N_r\}$
satisfying $dH^N_s/ds=2NV^N_s$ for almost all $s$, with $V^N_s\in\{-1,1\}$. The other identities which we will state in this proof are true a.s. In these identities we exclude the trajectories of $H^N$ which have a local maximum or minimum at the level $t$. This implies  that the two processes
$s\to V^N_s$ and $s\to {\bf1}_{\{H^N_s>t\}}$ do not jump at the same time. Hence from
\begin{align*}
{\bf1}_{\{H^N_s>t\}}&=\sum_{0<r<s}{\bf1}_{\{H^N_r=t\}}V^N_r-{\bf1}_{\{V^N_s=-1\}}{\bf1}_{\{H^N_s=t\}}\\
&=\frac{\sigma^2N}{2}\int_0^s V^N_rdL^N_r(t),
\end{align*}
we deduce by differentiating the product that
$$
V^N_s{\bf1}_{\{H^N_s>t\}}=\frac{\sigma^2 N}{2}\int_0^s(V^N_r)^2dL^N_r(t)
+\int_0^s{\bf1}_{\{H^N_r>t\}}dV^N_r.$$
Since $(V^N_r)^2=1$, this is the second identity  in the lemma.
\epf
\begin{lemma}\label{contLT}
Denote by $L_s(t)$ the local time at level $t$ up to time $s$ of $H$. Then with probability one $(s,t)\mapsto L_s(t)$ is continuous from $\mathbb{R}_{+}\times \mathbb{R}_{+}$ into $\mathbb{R}$.
\end{lemma}
\bpf
This is Theorem VI.1.7 page $225$ of Revuz, Yor \cite{ReYo}.
\epf

\begin{proposition}\label{convL(s)}
 For each $d\ge1$, $0\le t_1<t_2<\cdots<t_d$,
\begin{align*}
\big\{\big(H^N_s, L_{s}^N(t_1), L_{s}^N(t_2),..,L_{s}^N(t_d) \big), s \ge 0\big\} &\Rightarrow  \big\{\big( H_s,L_{s}(t_1), L_{s}(t_2),..,L_{s}(t_d) \big), s \ge 0\big\}\\ &\mbox{ in } C([0,\infty ))\times (D([0,\infty )))^d.
\end{align*}
\end{proposition}
\bpf
We prove the result in case $d=1$ only, the proof of the general case being very similar. 
From \eqref{HnVn} and Lemma \ref{naive} we deduce that for any $t\ge0$, a. s.
\begin{equation}\label{TanakaN}
L_s^N(t)= 2(H_s^N-t)^{+}+ \frac2{N\sigma ^2} V_s^N \mathbf{1}_{\{H_s^N> t\}}- 2 \int _0^s \mathbf{1}_{\{H_r^N> t\}} (dM_r^{1,N}-dM_r^{2,N}),\ s\ge0.
\end{equation}
Let 
\begin{equation*}
U_s^N = \int _0^s \mathbf{1}_{\{H_r^N> t\}} (dM_r^{1,N}-dM_r^{2,N}).
\end{equation*}
By Proposition \ref{sc} we have that $\{U^N\}_{N\geq 1}$ is tight in $D([0,\infty ))$. Moreover,
\begin{align}
\langle M^{1,N}- M^{2,N} \rangle_s &= \frac 4{\sigma ^2} s\label{M1-M2}\\
\langle U^N \rangle_s = \langle U^N,M^{1,N}- M^{2,N} \rangle_s &= \frac 4{\sigma ^2} \int_0^s \mathbf{1}_{\{H_r^N> t\}} dr.\label{covar}
\end{align}
From the occupation times formula
\begin{equation*}
\int_0^s \mathbf{1}_{\{H_r= t\}} dr = \frac {\sigma ^2}4 \int_0^{\infty } \mathbf{1}_{\{r= t\}} L_s(r) dr= 0 \quad a.s.
\end{equation*}
Then by Lemma \ref{convloctim} from the Appendix we deduce that along an appropriate sequence
\begin{equation}\label{convUN}
\Big \{ \int_0^s \mathbf{1}_{\{H_r^N> t\}} dr, s\geq 0 \Big\} \Rightarrow \Big \{ \int_0^s \mathbf{1}_{\{H_r> t\}} dr, s\geq 0 \Big\}
\end{equation}
From \eqref{M1-M2}, \eqref{covar} and \eqref{convUN}, we have again along an appropriate subsequence
\begin{equation*}
(U_s^N, M_s^{1,N}- M_s^{2,N}) \Rightarrow  \big( \frac2{\sigma } \int _0^s \mathbf{1}_{\{H_r> t\}} dB_r, \frac2{\sigma } B_s \big) \quad \mbox{in} \quad \big( D([0,\infty ))\big)^2.
\end{equation*}
Moreover, arguments similar to that used in the proof of Proposition \ref{convcombin}
establish that
 $$(H_s^N, U_s^N) \Rightarrow \big( H_s, \frac2{\sigma } \int _0^s \mathbf{1}_{\{H_r> t\}} dB_r \big)\quad \mbox{ in }\quad \big(D([0,\infty ))\big)^2.$$
Now from any subsequence, we can extract a subsequence along which we can take the weak limit in 
\eqref{TanakaN}. But Tanaka's formula gives us the identity
\begin{equation*}
L_s(t)= 2(H_s-t)^{+}- \frac4{\sigma } \int _0^s \mathbf{1}_{\{H_r> t\}} dB_r,
\end{equation*}
which characterizes the limit of $L^N$ as the local time of $H$. Since the law of $H$ is uniquely characterized,
the whole sequence converges.
\epf
 
\begin{proposition}\label{tightL(t)}
For each $s\geq 0$ fixed, $\{L_{s}^N(t), t\geq 0\}_{N\geq 1}$ is tight in $D([0,\infty ))$.
\end{proposition}
\bpf
We have
\begin{align}
L_{s}^N(t)&= 2(H_{s}^N-t)^{+}+ \frac2{N\sigma ^2} V_{s}^N \mathbf{1}_{\{H_{s }^N> t\}}+ \frac4{N \sigma ^2} \int _0^{s} V_{r-}^N \mathbf{1}_{\{H_r^N> t\}} dM_r^N\label{LN}\\
&= K_t^N + G_t^N,\nonumber
\end{align}
where
\begin{align*}
K_t^N&= 2(H_{s}^N-t)^{+}+ \frac2{N\sigma ^2} V_{s}^N \mathbf{1}_{\{H_{s}^N> t\}},\\
G_t^N&= \frac4{N \sigma ^2} \int _0^{s} V_{r-}^N \mathbf{1}_{\{H_r^N> t\}} dM_r^N.
\end{align*}
From
\begin{align*}
&\{K^N_0=2H^N_s+\frac{2}{N\sigma^2}V^N_s{\bf1}_{\{H^N_s>0\}},\ N\ge1\}\quad\text{is tight and}\\
&\limsup_{N\to\infty}|K^N_t-K^N_{t'}|\le2|t-t'|,
\end{align*}
it follows from Theorem 15.1 in \cite{Bi} that 
the sequence $\{K_{\cdot}^N\}_{N\geq 1}$ is tight, and any limit of a converging subsequence
is a. s. continuous.

We next show that the sequence $\{G_{\cdot}^N\}_{N\geq 1}$ satisfies the conditions of Proposition \ref{tightEK}.

Condition $(1)$ follows easily from the fact that $\mathbb{E} \big( |G_t^N|^2 \big) \leq 16s/ \sigma ^2$. In order to verify condition $(2)$, we will show that for any $T>0$, there exists $C>0$ such that for any $0<t<T, \varepsilon >0$,
\begin{equation*}
\mathbb{E} \big[ (G_{t+\varepsilon }^N- G_t^N)^2 (G_t^N- G_{t-\varepsilon }^N)^2 \big] \leq C(\varepsilon ^{3/2}+ \varepsilon ^2).
\end{equation*}
In order to simplify the notations below we let
\begin{align*}
\varphi _r^N&:= V_{r-}^N \mathbf{1}_{\{t-\varepsilon <H_r^N\leq  t\}},\\
\psi _r^N&:= V_{r-}^N \mathbf{1}_{\{t <H_r^N\leq  t+\varepsilon \}}.
\end{align*}
An essential property, which will be crucial below, is that $\varphi _r^N \psi _r^N=0$. Also $(\varphi _r^N)^2=|\varphi _r^N|$, and similarly for $\psi ^N$, since those functions take their values in the set $\{-1,0,1\}$. The quantity we want to compute equals up to a fixed multiplicative constant
\begin{equation*}
N^{-4}\mathbb{E} \Big[ (\int_0^s \varphi _r^N dM_r^N)^2 (\int_0^s \psi _r^N dM_r^N)^2 \Big].
\end{equation*}
We note that we have the identity
\begin{equation*}
(\int_0^s \varphi _r^N dM_r^N)^2= 2 \int_0^s \int_0^{r-} \varphi _u^N dM_u^N \varphi _r^N dM_r^N+ \int_0^s |\varphi_r^N| dM_r^N+ \sigma ^2N^2 \int_0^s |\varphi _r^N| dr,
\end{equation*}
and similarly with $\varphi ^N$ replace by $\psi ^N$. Because $\varphi _r^N \psi _r^N=0$, the expectation of the product of
\begin{equation*}
\int_0^s \int_0^{r-} \varphi _u^N dM_u^N \varphi _r^N dM_r^N \quad or \quad \int_0^s |\varphi_r^N| dM_r^N
\end{equation*}
with
\begin{equation*}
\int_0^s \int_0^{r-} \psi _u^N dM_u^N \psi _r^N dM_r^N \quad or \quad \int_0^s |\psi _r^N| dM_r^N
\end{equation*}
vanishes. We only need to estimate the expectations
\begin{align*}
&\mathbb{E} \Big( \int_0^s \int_0^{r-} \varphi _u^N dM_u^N \varphi _r^N dM_r^N \int_0^s |\psi _r^N| dr \Big), \mathbb{E} \Big( \int_0^s |\varphi_r^N| dM_r^N \int_0^s |\psi _r^N| dr \Big),\\
&and \quad \mathbb{E} \Big( \int_0^s |\varphi _r^N| dr \int_0^s |\psi _r^N| dr \Big),
\end{align*}
together with similar quantities with $\varphi ^N$ and $\psi ^N$ interchanged. The estimates of the first two expectations are very similar. We estimate the second one as follows, using the Cauchy-Schwarz inequality, and Lemma \ref{estimtt'} below~:
\begin{align*}
\mathbb{E} \Big( \int_0^s |\varphi_r^N| dM_r^N \int_0^s |\psi _r^N| dr \Big) &\leq \sqrt{\mathbb{E} \int_0^s |\varphi _r^N| d \langle M^N \rangle _r} \sqrt{\mathbb{E} \Big[ \Big( \int_0^s |\psi _r^N| dr \Big)^2 \Big]}\\
& \leq C N \varepsilon ^{3/2}
\end{align*}
Finally, again from Lemma \ref{estimtt'},
\begin{equation*}
\mathbb{E} \Big( \int_0^s |\varphi _r^N| dr \int_0^s |\psi _r^N| dr \Big) \leq C \varepsilon ^2.
\end{equation*}
The first quantity should be multiplied by $N^2$, and the second by $N^4$, and then both should be divided by $N^4$. The proposition now follows from Proposition \ref{tightcombined}.
\epf

\begin{lemma}\label{estimtt'}
Let $s,\epsilon ,T>0$. Then there exists a constant $C$ such that for all $N\geq 1$ and $0<t,t'< T$,
\begin{align*}
&\mathbb{E}\Big( \int_0^{s} \mathbf{1}_{\{t-\varepsilon <H_r^N \leq  t\}} dr \Big) \leq C \varepsilon, \\
&\mathbb{E}\Big( \int_0^{s} \mathbf{1}_{\{t-\varepsilon <H_r^N \leq  t\}} dr \int_0^{s} \mathbf{1}_{\{t'-\varepsilon <H_r^N \leq  t'\}} dr \Big) \leq C \varepsilon^2.
\end{align*}
\end{lemma}
\bpf
We will prove the second inequality, the first one follows from the second one with $t=t'$ and the Cauchy-Schwarz inequality.

For $s,t>0$ define $F_s^N(t):= \int_0^s \mathbf{1}_{\{0\leq H_r^N \leq  t\}} dr$. It follows readily from the definition of $L^N$ that
\begin{equation*}
\frac{\partial F_s^N}{\partial t} (t) = \frac{\sigma ^2}{4} L_s^N(t).
\end{equation*}
Hence
\begin{align*}
\mathbb{E}\Big( \int_0^{s} \mathbf{1}_{\{t-\varepsilon <H_r^N \leq  t\}} dr \int_0^{s} \mathbf{1}_{\{t'-\varepsilon <H_r^N \leq  t'\}} dr \Big)&= \frac{\sigma ^4}{16}\mathbb{E}\Big( \int_{t-\varepsilon }^t L_s^N(r) dr \int_{t'-\varepsilon }^{t'} L_s^N(u) du \Big)\\
&= \frac{\sigma ^4}{16}\mathbb{E}\Big( \int_{t-\varepsilon }^t  \int_{t'-\varepsilon }^{t'} L_s^N(r) L_s^N(u) dr du \Big)\\
&= \frac{\sigma ^4}{16} \int_{t-\varepsilon }^t  \int_{t'-\varepsilon }^{t'} \mathbb{E}\big( L_s^N(r) L_s^N(u) \big) dr du\\
&\leq \frac{\sigma ^4}{16} \varepsilon ^2 \sup_{0\leq r,u\leq T} \mathbb{E}\big( L_s^N(r) L_s^N(u) \big) \\
&=\frac{\sigma ^4}{16} \varepsilon ^2 \sup_{0\leq r\leq T} \mathbb{E}\big( (L_s^N(r))^2  \big).
\end{align*}
On the other hand, since
by It\^o's formula there exists a martingale $\bar {M}_s^N$ such that
\begin{equation*}
(H_s^N)^2+ \frac{2}{N \sigma ^2}H_s^N V_s^N= \frac4{\sigma ^2}s+ \bar {M}_s^N.
\end{equation*}
we conclude that 
\begin{equation*}
\sup_{N\geq 1} \mathbb{E} \big( (H_s^N)^2 \big) <\infty .
\end{equation*}
The second inequality now follows from \eqref{LN}.
\epf
\begin{proposition}\label{convL(t)}
 For all $d>1$, $0\le s_1<s_2<\cdots<s_d$,
\begin{equation*}
(H^N,L_{s_1}^N,\ldots,L_{s_d}^N)\Rightarrow (H,L_{s_1},\ldots,L_{s_d}) \quad in \quad
C([0,\infty))\times (D([0,\infty )))^d.
\end{equation*}
\end{proposition}
\bpf  We prove the result in the case $d=1$ only, the proof in the general case being very similar.
From
Proposition \ref{convL(s)} there follows in particular that 
for all $k\geq 1, 0\leq t_1<t_2<\cdots<t_k$, we have
\begin{equation*}
\big( H^N,L_{s}^N(t_1), L_{s}^N(t_2),..,L_{s}^N(t_k) \big) \Rightarrow  \big( H,L_{s}(t_1), L_{s}(t_2),..,L_{s}(t_k) \big)
\end{equation*}
That is, $\{L_{s}^N\}$ converges in finite-dimensional distributions to $\{L_{s}\}$, jointly with $H^N$. By Proposition \ref{tightL(t)}, $\{L_{s }^N(t), t\geq 0\}_{N\geq 1}$ is tight. The result follows.
\epf

We are now prepared to complete the \\
{\sc Proof of Theorem \ref{conHL}}: The main task is to combine the assertions of Propositions \ref {convL(s)} and \ref{convL(t)}, which means to turn the ``partial'' convergences asserted for $L^N$ in these propositions into a convergence that is joint in $s$ and $t$. We will also combine this result with Proposition \ref{convcombin} in order to get joint convergence of all our processes. To facilitate the reading, we will divide the proof  into several steps.

\bigskip 

\noindent
{\sc Step 1.}
Let $\{s_n, n\geq 1\}$ denote a countable dense subset of $\mathbb{R}_{+}$. Our first claim is that for all $n \in \mathbb N$,
\begin{align}\label{s1sn}
(H^N,M^{1,N},M^{2,N}, L_{s_1}^N,\ldots, L_{s_n}^N,S^N_x)\Rightarrow (H,\frac{2}{\sigma}B^1,\frac{2}{\sigma}B^2, L_{s_1},\ldots, L_{s_n},S_x)\\ \notag  \mbox{ in }  C(\mathbb{R}_{+}) \times D(\mathbb{R}_{+})^{n+2}\times \R_+.
\end{align}
To make the core of the argument  clear, let us write just for the moment
$$Y^N:=(M^{1,N},M^{2,N},S^N_x), \, Y:=(\frac{2}{\sigma}B^1,\frac{2}{\sigma}B^2,S_x), \, \Lambda^N:= (L_{s_1}^N,\ldots, L_{s_n}^N), \, \Lambda := (L_{s_1},\ldots, L_{s_n}).$$
Then \eqref{s1sn} translates into
\begin{equation}\label{abstr}
(H^N,Y^N, \Lambda^N)\Rightarrow (H,Y,\Lambda).
\end{equation}
By Proposition \ref{convcombin}, $(H^N,Y^N)\Rightarrow (H,Y)$, and by Proposition \ref{convL(t)}, $(H^N,\Lambda^N)\Rightarrow (H,\Lambda)$. Because in our situation $\Lambda$ is a.s. a function of $H$, these two convergences imply \eqref{abstr}. (More generally, this implication would be true if $Y$ and $\Lambda$ would be conditionally independent given $H$.)

\bigskip

\noindent
{\sc Step 2.} Now having established \eqref{s1sn}, it follows from a well known theorem due to Skorohod that all the processes appearing there can be constructed on a joint probability space, such that there exists an event $\mathcal{N}$ with $\mathbb{P}( \mathcal{N})=0$ and for all $\omega \notin \mathcal{N}$, 
\begin{equation}\label{convas1}
S^N_x(\omega)\to S_x(\omega),
\end{equation}
\begin{equation}\label{convas2}
(H^N_s(\omega),M^{1,N}_s(\omega ),M^{2,N}_s(\omega ))\rightarrow (H_s(\omega),\frac{2}{\sigma}B^1_s(\omega ), \frac{2}{\sigma}B^2_s(\omega )) \quad \textrm{locally uniformly in }s\geq 0,
\end{equation}
and for all $n\geq 1$,
\begin{equation}\label{convLunif}
L_{s_n}^N(t)(\omega )\rightarrow L_{s_n}(t)(\omega ) \quad \textrm{locally uniformly in }t\geq 0,
\end{equation}
as $N\rightarrow \infty $. Here we have made use of Lemma \ref{contLT}, which allows us to assume that $(s,t)\mapsto L_s(t)(\omega )$ is continuous from $\mathbb{R}_{+}\times \mathbb{R}_{+}$ into $\mathbb{R}$ for all $\omega \notin \mathcal{N}$, possibly at the price of enlarging the null set $\mathcal{N}$, and of Lemma \ref{CdansD} from the Appendix.

\bigskip

\noindent{\sc Step 3.} We claim that in the situation described in the previous step one even has
 for all $C,T>0, \omega \notin \mathcal{N}$,
\begin{equation}\label{locunif}
\sup_{0\leq s\leq C, \, 0\leq t\leq T} |L_s^N(t,\omega )-L_s(t,\omega )| \rightarrow 0,
\end{equation}
as $N\rightarrow \infty $. In other words, in Skorokhod's construction there is a.s. convergence of $L_s^N(t)$ to $L_s(t)$, locally uniformly in $s$ and $t$.  To prove \eqref{locunif}, we will make use of the fact that for any $\omega \notin \mathcal{N}$, and all $N,t$, the mapping $s\mapsto L_s^N(t)(\omega )$ is increasing and the mapping $s\mapsto L_s(t)(\omega )$ is continuous and increasing.
  Moreover, since the mapping $(s,t)\mapsto L_s(t,\omega )$ is continuous from the compact set $[0,C]\times [0,T]$ into $\mathbb{R}_{+}$, for any $\varepsilon >0$, there exists $\delta >0$ such that $0\leq s<s'\leq C, 0\leq t\leq T$ and $s'-s\leq \delta $ implies that
\begin{equation*}
L_{s'}(t,\omega )-L_s(t,\omega ) \leq \varepsilon .
\end{equation*}
Hence there exists $k\geq 1$ and $0=:s_0<r_1<\cdots<r_k:=C$ such that $\{r_i, 0\leq i < k\} \subset \{s_n, n\geq 1\}$ and moreover, $ r_i-r_{i-1}\leq \delta $  for all $1\leq i\leq k$. We have
\begin{equation*}
\sup_{0\leq s\leq C, 0\leq t\leq T} |L_s^N(t,\omega )-L_s(t,\omega )|\leq \sup_{1\leq i\leq k} [A_{N,i}+B_{N,i}]
\end{equation*}
where
\begin{align*}
&A_{N,i}= \sup_{r_{i-1}\leq s\leq r_i, 0\leq t\leq T} (L_s^N(t,\omega )-L_s(t,\omega ))^{+} \\
&B_{N,i}= \sup_{r_{i-1}\leq s\leq r_i, 0\leq t\leq T} (L_s^N(t,\omega )-L_s(t,\omega ))^{-}.
\end{align*}
For $r_{i-1}\leq s\leq r_i$, 
\begin{align*}
(L_s^N(t,\omega )-L_s(t,\omega ))^{+} &\leq (L_{r_i}^N(t,\omega )-L_s(t,\omega ))^{+} \\
& \leq (L_{r_i}^N(t,\omega )-L_{r_i}(t,\omega ))^{+} + \varepsilon , \\
(L_s^N(t,\omega )-L_s(t,\omega ))^{-} &\leq (L_{r_{i-1}}^N(t,\omega )-L_s(t,\omega ))^{-} \\
& \leq (L_{r_{i-1}}^N(t,\omega )-L_{r_{i-1}}(t,\omega ))^{-} + \varepsilon .
\end{align*}
Finally, 
\begin{equation*}
\sup_{0\leq s\leq C, 0\leq t\leq T} |L_s^N(t,\omega )-L_s(t,\omega )|\leq 2 \sup_{0\leq i\leq k} \sup_{0\leq t\leq T} |L_{r_i}^N(t,\omega )-L_{r_i}(t,\omega )| +2\varepsilon ,
\end{equation*}
while from \eqref{convLunif},
\begin{equation*}
\limsup_{N\rightarrow \infty } \sup_{0\leq s\leq C, 0\leq t\leq T} |L_s^N(t,\omega )-L_s(t,\omega )|\leq 2\varepsilon .
\end{equation*}
This implies \eqref{locunif}, since $\varepsilon >0$ is arbitrary. The assertion of Theorem \ref{conHL} is now immediate by combining \eqref{convas1}, \eqref{convas2} and \eqref{locunif}.
\epf

\section{Change of measure and proof of Theorem \ref{main}}\label{secgirs}

As in the previous section, let, for fixed $N\in \mathbb N$, $H^N$ be a process that follows the dynamics described in Proposition \ref{dynHN} for $\theta = \gamma =0$. We denote the underlying probability measure by $\mathbb P$, and the filtration by $\mathcal F = (\mathcal F_s)$.   Our first aim is to construct, by a Girsanov reweighting of the restrictions $\mathbb P|_{\mathcal F_s}$, a measure $\tilde {\mathbb P}^N$ under which $H^N$ follows the dynamics from Proposition \ref{dynHN} for a prescribed $\theta \ge 0$ and $\gamma >0$.

Here, a crucial role is played by the point process $P^N$ of the successive local maxima and minima of $H^N$, excluding the minima at height $0$.  Under $\mathbb P$, this is a Poisson process with intensity $\sigma^2 N^2$. 
More precisely, the process $Q^{1,N}$ which counts the successive local minima of $H^N$ (except those at height $0$) is a point process with predictable intensity $\lambda_s^{1,N}:= N^2\sigma^2 \mathbf 1_{\{V^N_{s-}=-1\}}$, and the process $Q^{2,N}$ which counts the successive local maxima of $H^N$ is a point process with predictable intensity $\lambda_s^{2,N}:= N^2\sigma^2 \mathbf 1_{\{V^N_{s-}=+1\}}$.
(Recall that the process $V^N$ is the (c\`adl\`ag) sign of the slope of $H^N$.)

For the rest of this section we fix $\theta \ge 0$ and $\gamma >0$. In view of Proposition \ref{dynHN} we want to change the rate $\lambda_s^{1,N}$ to $\lambda_s^{1,N}(1+\frac{2\theta}{N\sigma^2})$ and the rate $\lambda_s^{2,N}$ to $\lambda_s^{2,N}(1+\frac{4\gamma \Lambda^N_s(H_s)}{N\sigma^2})$.

As in Section \ref{allzero} we will use the process $M^N_s = P^N_s-N\sigma^2 s, \, s\ge 0$, which is a martingale under $\mathbb P$. Taking the route designed by Proposition \ref{Gir-PP} in the Appendix, we consider the local martingales
$$X_s^{N,1}:= \int_0^s\frac{2\theta }{N\sigma ^2} \mathbf{1}_{\{V_{r-}^N=-1\}}dM_r^N, \quad X_s^{N,2}:=  \int_0^s\frac{\gamma L_r^N(H_r^N)}{N}\mathbf{1}_{\{V_{r-}^N=1\}}dM_r^N, \quad X^N := X^{N,1}+X^{N,2}. $$
Let $Y^N:= \mathcal E(X^N)$ denote the Dol\'eans exponential 
of $X^N$.  Proposition \ref{expdol} in the Appendix recalls this concept  and the fact that  $Y^N$ is the solution of
\begin{equation}\label{eqY}
Y_s^N = 1+ \int_0^s Y_{r-}^N \left(\frac{2\theta }{N\sigma ^2} \mathbf{1}_{\{V_{r-}^N=-1\}}+\frac{\gamma L_r^N(H_r^N)}{N}\mathbf{1}_{\{V_{r-}^N=1\}}\right)dM_r^N,  \quad s\ge 0.
\end{equation}
We will show that $Y^N$ is a martingale, which from Proposition \ref{Gir-PP} will directly render the required change of measure. 
\begin{proposition} \label{PNtilde}
$Y^N$ is a $(\mathcal F, \mathbb P)$-martingale. 
\end{proposition}
\bpf Under $\mathbb P$, $Y^N$ is a positive super--martingale and a local martingale. 
It is a martingale if and only if 
\begin{equation}\label{Ymart}
\mathbb{E}[Y_s^N] =Y^N_0= 1,
\end{equation}
 which we will show.
A key idea is to work along the excursions of $H^N$, that is, along the sequence of stopping times $\tau_a^{N,s}:= S_{a/N}^N\wedge s$, $a=0,1,2,\ldots$ Since  $N$ and $s$ are fixed, we will suppress the superscripts $N$ and $s$ for brevity and write $\tau_a$ instead of $\tau_a^{N,s}$.  

\bigskip  

\noindent{\sc Step 1.} We want to show that
\begin{align}\label{Ytauamart}
\mathbb E[Y^N_{\tau_a}] = 1, \quad a=1,2,\ldots
\end{align}
Between $\tau_{a-1}$ and $\tau_a$, the solution of \eqref{eqY} is bounded above by the solution
of the same equation with $M^N$ replaced by $P^N$, which takes the form $d\tilde{Y}^N_r=\tilde{Y}^N_{r-}a^N_rdP^N_r$.
If we denote by $\{T_k,\ k\ge1\}$ the successive jump times of the Poisson process $P^N$, we have or each $k$
$\tilde{Y}^N_{T_k}=\tilde{Y}^N_{T_k-}(1+a^N_{T_k})$. Consequently
 for all $\tau_{a-1}\le r\le\tau_a$,
\begin{equation*}
\frac{Y^N_r}{Y^N_{\tau_{a-1}}} \leq \prod _{k\geq 1: \tau_{a-1}\le T_k\leq \tau_a}\big(1+ \frac{2\theta }{N\sigma ^2} \mathbf{1}_{\{V_{T_k-}^N=-1\}}+\frac{\gamma L_{T_k}^N(H_{T_k}^N)}{N}\mathbf{1}_{\{V_{T_k-}^N=1\}}\big).
\end{equation*}
Within the excursion of $H^N$ between the times $\tau_{a-1}$ and $\tau_a$, these jump times coincide with the times of the local maxima and minima of $H^N$ in the time interval $(\tau_{a-1}, \tau_a)$.  Since for $a>1$ there are reflections of $X^N$ at 0 in the time interval $(0, \tau_{a-1})$, the parity of those $k$ for which $V_{T_k-}^N=1$, $\tau_{a-1}\le T_k\leq \tau_a$, depends on $a$. However,   
noting that 
\begin{equation}\label{boundLH}
L_{T_k}^N(H_{T_k}^N)\leq \frac{4}{N \sigma ^2} k,
\end{equation}
 we infer the existence of a constant $c>0$ such that 
\begin{equation}\label{estRatio}
\frac{Y^N_r}{Y^N_{\tau_{a-1}}}  \leq   c^{P^N_{\tau_a}} \times  (P^N_{\tau_a}+1)!!\, , \quad
\tau_{a-1}\le r\le\tau_a,
\end{equation}
where we define for $k\in \mathbb{N}, k!!= 1\cdot 3\cdot 5\cdots k$ if $k$ is odd, and $k!!= 1\cdot 3\cdot 5\cdots(k-1)$ if $k$ is even. 

Now
\begin{equation}\label{estimY}
Y^N_{\tau_a}=Y^N_{\tau_{a-1}}\left(1+\int_{\tau_{a-1}}^{\tau_a}\frac{Y^N_{r-}}{Y^N_{\tau_{a-1}}} [p(r)+L^N_{r-}(H^N_r)q(r)]dM^N_r\right),
\end{equation}
where $0\le p(r)\le{\frac2\theta}{N^2}$, $0\le q(r)\le\frac{\gamma}{N}$, $p$ and $q$ are predictable.
 The claimed equalities $\mathbb{E}[Y^N_{\tau_a}]=1$, $a=1,2,\ldots$,  follow by induction on $a$
from \eqref{estimY}, provided the process
$$\mathcal{M}^N_s=\int_0^s{\bf1}_{]\tau_{a-1},\tau_a]}(r)\frac{Y^N_{r-}}{Y^N_{\tau_{a-1}}} [p(r)+L^N_{r-}(H^N_r)q(r)]dM^N_r, \quad s\ge 0,$$ is a martingale. From Theorem T$8$ in Br\'emaud \cite{PB} page $27$, this is a consequence of the fact that
$$\E\int_0^s{\bf1}_{]\tau_{a-1},\tau_a]}(r)\frac{Y^N_{r-}}{Y^N_{\tau_{a-1}}} [p(r)+L^N_{r-}(H^N_r)q(r)]dr<\infty.$$
In order to verify the latter inequality, we compute
\begin{align*}
\E\left[\int_{\tau_{a-1}}^{\tau_a}\frac{Y^N_{r-}}{Y^N_{\tau_{a-1}}}[p(r)+L^N_{r-}(H^N_r)q(r)] dr\right]
&\le C_N s \E\left[ c^{P^N_{\tau_a}}\times(P^N_{\tau_a}+1)!!(1+P^N_{\tau_a})\right]\\
&\le C_N s \E\left[ c^{P^N_s}\times(P^N_s+1)!!(1+P^N_s)\right]\\
&\le C_N s C_{N,s} ,
\end{align*}
where we have used \eqref{estRatio}, \eqref{boundLH} and $\tau_a\le s$, and
where $C_N$ and $C_{N,s}$ are constants which depend only on $N$ and $(N,s)$, respectively.
The fact that $C_{N,s}<\infty$ follows from
\begin{equation*}
\mathbb{E}[c_2^{P^N_s}(P^N_s+1)!! P^N_s]= \exp\{-N^2\sigma ^2 s\} \sum_{k=0}^{\infty } c_2^k(k+1)!! k \frac{(N^2\sigma ^2 s)^k}{k!}
\end{equation*}
Since 
\begin{align*}
\frac{(k+1)!!k}{k!}= \frac{k(k+1)}{2.4...\big(2[\frac{k}2]\big)}&=\frac{k(k+1)}{2^{[\frac{k}2]}}. \frac1{[\frac{k}2]!}\\
&< \quad \frac1{[\frac{k}2]!} \quad \quad \forall k\geq 20,
\end{align*}
we deduce that $\mathbb{E}[c_2^{P^N_s}(P^N_s+1)!! P^N_s]<\infty $. This completes the proof of \eqref{Ytauamart}.

\bigskip

\noindent{\sc Step 2.}
 We can now define a consistent family of probability measures $\tilde {\mathbb{P}}^{N,s,a}$ on $\mathcal{F}_{\tau_a}$, $a=1,2,\ldots$ by putting
\begin{equation*}
\frac{d\tilde {\mathbb{P}}^{N,s,a}}{d\mathbb{P}\mid _{\mathcal{F}_{\tau_a}}}= Y_{\tau_a}^N, \quad a\in \mathbb N.
\end{equation*}
We write $ \tilde{\mathbb{P}}^{N,s}$ for the probability measure on the $\sigma$-field generated by union of the $\sigma$-fields $\mathcal F_{\tau_a}$, $a=1,2,\ldots$, whose restriction to $\mathcal{F}_{\tau_a}$ is $\tilde {\mathbb{P}}^{N,s,a}$ for all $a=1,2,\ldots$,  and put
$$A:= \inf\{a \in \mathbb  N:\tau_a =s\}.$$
We will now show that 
\begin{enumerate}
\item{(i)} $A < \infty$ \,  $\tilde{\mathbb{P}}^{N,s}$- a.s. (and consequently $\tau_A = s$ \,  $\tilde{\mathbb{P}}^{N,s}$- a.s.),   
\item{(ii)} under $\tilde {\mathbb{P}}^{N,s}$, $(H^N_r)_{0\le r \le \tau_A}=(H^N_r)_{0\le r \le s}$ is a stochastic process following the dynamics specified in Proposition \ref{dynHN}. 
\end{enumerate}
Indeed, applying Girsanov's theorem (Proposition  \ref{Gir-PP} in the Appendix) to the $2$-variate point process 
\begin{align}\label{Q12}
(Q_r^{1,N},Q_r^{2,N})=\left( \int_0^r \mathbf{1}_{\{V_{u-}^N=-1\}} dP^N_u, \int_0^r \mathbf{1}_{\{V_{u-}^N=1\}} dP^N_u \right), \quad 0\le r \le \tau_a,
\end{align}
we have that under $\tilde{\mathbb{P}}^{N,s,a}$
\begin{align*}
&Q_r^{1,N}\quad \textrm{has intensity}\quad (N^2\sigma ^2+ 2\theta N)\mathbf{1}_{\{V_{r-}^N=-1\}}dr\\
&Q_r^{2,N}\quad \textrm{has intensity}\quad \sigma ^2[N^2+ \gamma N L_r^N(H_r^N)]\mathbf{1}_{\{V_{r-}^N=1\}}dr.
\end{align*}
Thus, for all $a\in \mathbb N$,  $(H^N_r)_{0\le r \le \tau_a}$ is,  under  $\tilde{\mathbb{P}}^{N,s,a}$, a stochastic process following the dynamics from Proposition \ref{dynHN} up to the stopping time $\tau_a$. Considering the sequence of excursions $(H^N_r)_{\tau_{a-1}\le r \le \tau_a}$, $a=1,2,\ldots$ under $ \tilde{\mathbb{P}}^{N,s}$, we infer from Lemma~\ref{extendable} the validity of the claims (i) and (ii).

\bigskip

\noindent{\sc Step 3.} We now prove \eqref{Ymart}. For this we observe that
$$\mathbb E [Y^N_s]= \sum_{a\ge 1} \mathbb E [Y^N_s; A = a] = \sum_{a\ge 1} \mathbb E [Y^N_{\tau_a}; A = a]=  \sum_{a\ge 1}\tilde{\mathbb{P}}^{N,s} (A=a) =  \tilde{\mathbb{P}}^{N,s} (A<\infty)=1.$$
\epf
\begin{corollary}\label{corPNtilde}
Let ${\tilde {\mathbb P}}^N$ be the probability measure on $\mathcal F$ whose restriction to $\mathcal F_s$, $s > 0$, has density $Y^N_s$ (given by \eqref{eqY}) with respect to $\mathbb P|_{\mathcal F_s}$. Then  under ${\tilde {\mathbb P}}^N$ the process $H^N$ follows the dynamics from Proposition \ref{dynHN} for the prescribed $\theta$ and $\gamma$.
\end{corollary}
\bpf This is immediate from Proposition \ref{PNtilde} and the discussion preceding it, combined with Proposition \ref{Gir-PP} in the Appendix applied to the process defined in \eqref{Q12}, now with $0\le r < \infty$.
\epf

Next we will analyze the behaviour of the Girsanov densities as $N\to \infty$. For this we use the two martingales $M^{1,N}$ and $M^{2,N}$ defined in \eqref{twomart}, and note that \eqref{eqY} can be rewritten as
$$Y_s^N =1+\int _0^s Y_{r-}^N \big\{\theta dM_r^{1,N}+\frac{\gamma\sigma^2  L_r^N(H_r^N)}{2} dM_r^{2,N} \big\}, \quad s\ge 0.$$
The two (pure jump) martingales $M^{1,N}$ and $M^{2,N}$ have jump sizes $2/(N\sigma^2)$, hence the random variable  under the expectation in formula \eqref{UCV} vanishes for suitably large $N$. Thus (see Definition \ref{defUCV} in the Appendix),  the sequences $\{M^{1,N}\}_{N\geq 1}$  and $\{M^{2,N}\}_{N\geq 1}$ have uniformly controlled variations, and because of Proposition \ref{th-good} (1) they are ``good''.  Hence
\begin{equation*}
X_\cdot^N \Rightarrow \int _0^\cdot \big\{\frac{\sqrt{2}\theta }{\sigma } dB_r^1+\frac{\sqrt{2}\gamma \sigma L_r(H_r)}{2} dB_r^2 \big\}:= X_\cdot.
\end{equation*}
Moreover, by Proposition \ref{th-good} (3), $\{X_s^N\}_{N\geq 1}$ is also a good sequence, hence by  Proposition \ref{th-good} (2)
 $$ Y^N = \mathcal{E}(X^N) \Rightarrow  \mathcal{E}(X)=: Y.$$ 
Combining these facts with Corollary \ref{convH-LS}, we deduce, again from Proposition \ref{th-good} (3), that
\begin{equation}\label{convH-LS-Y}
(H^N,L^N_{S^N_x},Y^N)\Rightarrow (H,L_{S_x},Y).
\end{equation}

Since $B^1$ and $B^2$ are mutually orthogonal, by Proposition \ref{expdol} we have
\begin{align*}
Y_s&= \mathcal{E}\Big( \int _0^{.} \big\{\frac{\sqrt{2}\theta }{\sigma } dB_r^1+\frac{\sqrt{2}\gamma \sigma L_r(H_r)}{2} dB_r^2 \big\}\Big)_s \\
&=\mathcal{E}\Big(\frac{\sqrt{2}\theta }{\sigma } B^1\Big)_s\, \,\mathcal{E}\Big(\int_0^{.}\frac{\sqrt{2}\gamma \sigma L_r(H_r)}{2} dB_r^2 \Big)_s \\
&=\exp\Big\{\frac{\sqrt{2}\theta }{\sigma } B_s^1+ \int_0^{s}\frac{\sqrt{2}\gamma \sigma L_r(H_r)}{2} dB_r^2-\int_0^s \big[\frac{\theta ^2}{\sigma ^2}+ \frac{\gamma ^2\sigma ^2}{4} L_r(H_r)^2\big] dr\Big\}
\end{align*}
Applying Lemma \ref{suffGirs} and Lemma \ref{Girstandard} again, we deduce that $Y$ is a martingale. In particular $\mathbb{E}[Y_s]=1$ for all $s\geq 0$.  Define the probability measure $\tilde {\mathbb{P}}$ by
\begin{equation*}
\frac{d\tilde {\mathbb{P}}\mid _{\mathcal{F}_s}}{d\mathbb{P}\mid _{\mathcal{F}_s}}= Y_s, \quad \forall s\geq 0,
\end{equation*}
then $H$, under $\tilde {\mathbb{P}}$, solves the SDE \eqref{reflecSDE} with $B_s$ there replaced by 
$$\tilde B_s := \frac 1{\sqrt 2} (B_s^1-  B_s^2) - \frac \theta \sigma s + \frac{\gamma \sigma} 2 \int_0^sL_r(H_r)dr,$$
 which is a standard Brownian motion under $\tilde {\mathbb{P}}$ due to Proposition \ref{Gir-BM}.
 
 The following general and elementary Lemma will allow us to conclude the required convergence under the transformed measures.
\begin{lemma}\label{weakConv-Girsanov}
Let $(\xi_N,\eta_N)$, $(\xi,\eta)$ be random pairs defined on a probability space $(\Omega,\mathcal{F},\mathbb{P})$,
with $\eta_N$, $\eta$ nonnegative scalar random variables, and $\xi_N$, $\xi$ taking values in some complete separable metric space $\mathcal{X}$.  Assume that $\E[\eta_N]=\E[\eta]=1$. Write $(\tilde{\xi}_N,\tilde{\eta}_N)$ 
for the random pair $(\xi_N,\eta_N)$ defined under the probability measure $\tilde{\mathbb P}^N$ which has density $\eta_N$ with respect to $\P$, and $(\tilde{\eta},\tilde{\xi})$ for the random pair $(\eta,\xi)$ defined under the probability measure $\tilde{\P}$ which has density $\eta$ with respect to $\P$. Then $(\tilde{\xi}_N,\tilde{\eta}_N)$ converges in distribution to
$(\tilde{\eta},\tilde{\xi})$, provided that $(\xi_N,\eta_N)$ converges in distribution to $(\xi,\eta)$.
\end{lemma}
\bpf
Due to the equality $\E[\eta_N]=\E[\eta]=1$ and a variant of Scheff\'e's theorem (see Thm. 16.12 in \cite{PBi}), the sequence ${\eta_N}$ is
uniformly integrable. Hence for all bounded continuous $F : \mathcal{X}\times\R_+\to\R$,
$$\E[F(\tilde{\xi}_N,\tilde{\eta}_N)]=\E[F(\xi_N,\eta_N)\eta_N]\to\E[F(\xi,\eta)\eta]=\E[F(\tilde{\xi},\tilde{\eta})].$$
\epf

Combining \eqref{convH-LS-Y} with Lemma \ref{weakConv-Girsanov} yields the
\begin{theorem}\label{th3}
Let $H^N$ be a stochastic process following the dynamics specified in Proposition \ref{dynHN}, and let $H$ be the unique weak solution of the SDE \eqref{reflecSDE}. We have
\begin{equation}\label{convHLS}
(H^N, L_{S_x^N}^N)\Rightarrow (H, L_{S_x})\quad \mbox {in} \quad C([0,\infty ])\times D([0,\infty ]),
\end{equation}
where $S_x^N$ and $S_x$ are defined in \eqref{defSN} and \eqref{Sx}.
\end{theorem}

We can now proceed with the

\noindent{\sc Completion of the proof of Theorem \ref{main}} :
Define $Z_t^{N,x}:= \frac{\sigma ^2}4 L_{S_x^N}^N(t)$. By Corollary~\ref{discreteRN},  $Z^{N,x}$ follows the dynamics \eqref{Zndyn}. From \eqref{convHLS}, $\frac{\sigma ^2}4 L_{S_x}$  is the limit in distribution of $Z^{N,x}$  as $N\rightarrow \infty $. Hence by Proposition \ref{convZ}, $t\mapsto \frac{\sigma ^2}4 L_{S_x}(t)$ is a weak solution of the SDE \eqref{fellog1}, which completes the proof of Theorem~\ref{main}. \epf


\begin{remark}
Theorem \ref{main} establishes a correspondence between the solution $H$ of the SDE
\eqref{reflecSDE} and the logistic Feller process, i.e.\ the solution of \eqref{fellog1}. This connection can be expressed in particular through the occupation times formula for $H$, which states that for any Borel measurable and positive valued function $f$,
$$\int_0^{S_x}f(H_s)ds=\int_0^\infty f(t)Z^x_tdt.$$
This formula in the particular case $f\equiv1$ states that 
$$S_x=\int_0^\infty Z^x_tdt.$$
The quantity on the right is the area under the trajectory $Z^x$. It is the limit of the 
properly scaled total branch length of the approximating forests $F^N$ defined in Section \ref{discretemass}. We now establish another identity concerning this same quantity, with the help of a time change introduced by Lambert in \cite{AL}. Consider the additive functional
$$A_t=\int_0^tZ^x_rdr,$$
and the associated time change
$$\alpha_t=\inf\{r>0,\ A_r>t\}.$$
As noted in \cite{AL}, the process $U^x_t:=Z^x_{\alpha_t}$ is an Ornstein--Uhlenbeck
process, solution of the SDE
$$dU^x_t=(\theta-\gamma U^x_t)dt+\sigma dB_t,\quad U^x_0=x.$$
Of course this identification is valid only for $0\le t\le\tau_x$, where
$\tau_x:=\inf\{t>0,\ U_t=0\}$. Let $T_x$ be the extinction time of the logistic Feller process
$Z^x_t$. We clearly have $\alpha_{\tau_x}=T_x$, and consequently
$$\tau_x=\int_0^\infty Z_r^xdr.$$
We have identified the time at which the local time at $0$ of the exploration process $H$ reaches
$x$ with the area under the logistic Feller trajectory starting from $x$, and with the time taken by the Ornstein--Uhlenbeck process $U^x$ to reach 0. The reader may notice that in the particular case $\gamma=0$, the identity $S_x=\tau_x$ is not a surprise, see also the discussion and the references in  \cite{PW} Section 6.
\end{remark}

\section{Appendix}
\subsection{Skorohod's topology and tightness in $D([0,\infty))$}
We denote by $D([0,\infty ))$ the space of functions from $[0,\infty )$ into $\mathbb{R}$ which are right continuous and have left limits at any $t>0$ (as usual such a function is called c\`adl\`ag). We briefly write $\mathbb{D}$ for the space of adapted, c\`adl\`ag stochastic processes. We shall always equip the space $D([0,\infty ))$ with the Skorohod topology, for the definition of which we refer the reader to Billingsley \cite{Bi} or Joffe, M\'etivier \cite{JoMe}. The next Lemma follows from 
considerations which can be found in \cite{Bi}, bottom of page 124.
\begin{lemma}\label{CdansD}
Suppose $\{x_n, n\geq 1\}\subset D([0,\infty ))$ and $x_n\rightarrow x$ in the Skorohod topology. 
\begin{description}
\item {(i)} If $x$ is continuous, then $x_n$ converges to $x$ locally uniformly.
\item {(ii)} If  each $x_n$ is continuous, then so is $x$, and $x_n$ converges to $x$ locally uniformly.
\end{description}
In particular, the space $C([0,\infty ))$ is closed in $D([0,\infty ))$ equipped with the Skorohod topology. 
\end{lemma}
The following two lemmata are used in the proofs of Propositions \ref{convL(s)} and \ref{convcombin}:
\begin{lemma}\label{convloctim}
Fix $t>0$. Let $x_n, x \in C([0,\infty )), n\geq 1$ be such that
\begin{description}
\item[\rm{1.}]$x_n\rightarrow x$ locally uniformly, as $n\to\infty$.
\item[\rm{2.}]for each $s>0$,
\begin{equation*}
\int_0^s \mathbf{1}_{\{x(r)= t\}} dr = 0.
\end{equation*}
\end{description}
Then 
\begin{equation*}
\int_0^s \mathbf{1}_{\{x_n(r)> t\}} dr \rightarrow \int_0^s \mathbf{1}_{\{x(r)> t\}} dr \quad \textrm{locally uniformly in $s\geq 0$.}
\end{equation*}
\end{lemma}
\bpf
We prove convergence for each $s>0$. The local uniformity is then easy. Given $\varepsilon >0$,  there exists $N_0$ such that
\begin{equation*}
\sup_{0\leq r\leq s} |x_n(r)-x(r)| < \varepsilon \quad \forall n\geq N_0.
\end{equation*}
Then for all $n\geq N_0$,
\begin{align*}
|\mathbf{1}_{\{x_n(r)> t\}} - \mathbf{1}_{\{x(r)> t\}}| &\leq \mathbf{1}_{\{t-\varepsilon < x(r)< t+\varepsilon \}}\\
\Big{|} \int_0^s \mathbf{1}_{\{x_n(r)> t\}} dr-\int_0^s \mathbf{1}_{\{x(r)> t\}} dr\Big{|} &\leq \int_0^s \mathbf{1}_{\{t-\varepsilon < x(r)< t+\varepsilon \}} dr.
\end{align*}
The result follows from
\begin{equation*}
\lim_{\varepsilon \rightarrow 0} \int_0^s \mathbf{1}_{\{t-\varepsilon < x(r)< t+\varepsilon \}} dr = \int_0^s \mathbf{1}_{\{x(r)= t\}} dr = 0.
\end{equation*}
\epf
\begin{lemma}\label{convStieltjes}
Let $x_n,y_n \in D([0,\infty )), n\geq 1$ and $x,y \in C([0,\infty ))$ be such that
\begin{description}
\item[\rm{1.}] for all $n \ge 1$, the function $t\to y_n(t)$ is increasing;
\item[\rm{2.}]$x_n\rightarrow x$ and $y_n\rightarrow y$, both locally uniformly.
\end{description}
Then $y$ is increasing and 
\begin{equation*}
\int_0^t x_n(s) dy_n(s) \rightarrow \int_0^t x(s) dy(s), \quad \textrm{locally uniformly in $t\geq 0$.}
\end{equation*}
\end{lemma}
\bpf
We prove convergence for each $t>0$. The local uniformity is then easy.
\begin{align*}
&\Big{|} \int_0^t x(s) dy(s)- \int_0^t x_n(s) dy_n(s)\Big{|} \\
&\leq \Big{|}\int_0^t [x(s)-x_n(s)]dy_n(s)\Big{|}+ \Big{|}\int_0^t x(s) [dy(s)- dy_n(s)]\Big{|}\\
&\leq \sup_{0\leq s\leq t} |x(s)- x_n(s)| y_n(t) + \int_0^t |x(s)- \xi _{\varepsilon }(s)| [dy(s)- dy_n(s)]+ \int_0^t |\xi _{\varepsilon }(s)| [dy(s)- dy_n(s)],
\end{align*}
where $\xi _{\varepsilon }$ is a step function which is such that $\sup_{0\leq s\leq t} |x(s)- \xi _{\varepsilon }(s)|\leq \varepsilon $. The first and last term of the above right hand side clearly tend to $0$ as $n\rightarrow \infty $. Then
\begin{align*}
\limsup_{n\rightarrow \infty } \Big{|} \int_0^t x(s) dy(s)- \int_0^t x_n(s) dy_n(s)\Big{|} &\leq \varepsilon \limsup_{n\rightarrow \infty } [y_n(t)+ y(t)]\\
&\leq 2y(t)\times \varepsilon .
\end{align*}
It remains to let $\varepsilon \rightarrow 0$.
\epf
We first state a tightness criterion, which is Theorem 13.5 from \cite{Bi}~:
\begin{proposition}\label{tightEK}
Let $\{X^n_t,\ t\ge0\}_{n\ge1}$ be a sequence of random elements of 
$D([0,\infty))$. A sufficient condition for $\{X^n\}$ to be tight is that the
two conditions (i) and (ii) be satisfied~:
\begin{description}
\item{(i)} For each $t \ge 0$, the sequence of random variables $\{X^n_t,\ n\ge1\}$ is tight in $\R$;
\item{(ii)} for each $T>0$, there exists $\beta, C>0$ and $\theta>1$ such that
$$\E\left(\left|X^n_{t+h}-X^n_t\right|^\beta\left|X^n_{t}-X^n_{t-h}\right|^\beta\right)\le C h^\theta,$$
for all $0\le t\le T$, $0\le h\le t$, $n\ge1$.
\end{description}
\end{proposition}

Note that convergence in $D([0,\infty))$ is not additive : $x_n\to x$ and $y_n\to y$ in 
$D([0,\infty))$ does not imply that $x_n+y_n\to x+y$ in $D([0,\infty))$. This is due to the fact that to the sequence $x_n$ is attached a sequence of time changes, and to the sequence $y_n$ is attached another sequence of time changes, such that the time changed $x_n$ and $y_n$ converge uniformly. But there may not exist a sequence of time changes which makes $x_n+y_n$ converge. If now $\{X^n_t,\ t\ge0\}_{n\ge1}$ and
$\{Y^n_t,\ t\ge0\}_{n\ge1}$ are two tight sequences of random elements of 
$D([0,\infty))$, we cannot conclude that $\{X^n_t+Y^n_t,\ t\ge0\}_{n\ge1}$ is tight.
However, if $x_n\to x$ and $y_n\to y$ in $D([0,\infty))$ and $x$ is continuous, then we deduce easily from Lemma 
\ref{CdansD} that $x_n+y_n\to x+y$ in $D([0,\infty))$.
It follows
 \begin{proposition}\label{tightcombined}
 If $\{X^n_t,\ t\ge0\}_{n\ge1}$ and
$\{Y^n_t,\ t\ge0\}_{n\ge1}$ are two tight sequences of random elements of 
$D([0,\infty))$ such that any limit of a weakly converging subsequence of
the sequence $\{X^n_t,\ t\ge0\}_{n\ge1}$  is a. s. continuous, then 
$\{X^n_t+Y^n_t,\ t\ge0\}_{n\ge1}$ is tight in $D([0,\infty))$.
\end{proposition}

Consider a sequence $\{X_t^n, t\geq 0\}_{n\geq 1}$ of one-dimensional semi--martingales, which is such that for each $n\geq 1$,
\begin{align*}
& X_t^n= X_0^n+ \int_0^t\varphi _n(X_s^n)ds+ M_t^n, \qquad t\geq 0; \\
& \langle M^n \rangle_t= \int_0^t\psi _n(X_s^n)ds, \qquad t\geq 0;
\end{align*}
where for each $n\geq 1, M^n$ is a locally square-integrable martingale, $\varphi _n$ and $\psi _n$ are Borel measurable functions from $\mathbb{R}$ into $\mathbb{R}$ and $\mathbb{R}_{+}$ respectively.

The following result is an easy consequence of Theorems 16.10 and 13.4 from \cite{Bi}.
\begin{proposition}\label{sc}
A sufficient condition for the above sequence $\{X_t^n, t\geq  0\}_{n\geq 1}$ of semi--martingales to be tight in $D([0,\infty ))$ is that both
\begin{equation}\label{tight0}
\textrm{the sequence of r.v.'s } \{X_0^n, n\geq  1\} \textrm{ is tight;}
\end{equation}
and for some $p>1$,
\begin{equation}\label{tightAC}
\forall T>0, \textrm{ the sequence of r.v.'s }\big\{\int_0^T[\mid \varphi _n(X_t^n)\mid +\psi _n(X_t^n)]^pdt, n\geq 1\big\}\textrm{ is tight.}
\end{equation}
Those conditions imply that both the bounded variation parts $\{V^n, n\geq 1\}$ and the martingale parts $\{M^n, n\geq 1\}$ are tight, and that the limit of any converging subsequence of $\{V^n\}$ is a.s. continuous.

If moreover, for any $T>0$, as $n\rightarrow \infty $,
\begin{center}
$\sup_{0\leq t\leq T}\mid M_t^n-M_{t-}^n\mid \rightarrow 0\quad$ in probability,
\end{center}
then any limit $X$ of a converging subsequence of the original sequence $\{X^n\}_{n\geq 1}$ is a.s. continuous.
\end{proposition}
\begin{remark}\label{remTight}
A sufficient condition for \eqref{tightAC}  is that for all $T>0$,
\begin{equation}\label{semiWeak}
\{\sup_{0\le t\le T}[|\varphi_n(X^n_t)|+\psi_n(X^n_t)],\ n\ge 1\}\text{ is tight}.
\end{equation}
\end{remark}
 \begin{remark}\label{tightnessCriteria}
 A sufficient condition for \eqref{tightAC}  is that for all $T>0$
 \begin{equation}\label{tightnessBound}
  \limsup_{n\ge1}\sup_{0\le t\le T}\E[\varphi^2_n(X^n_t)+\psi^2_n(X^n_t)]<\infty.
  \end{equation}
Indeed,  \eqref{tightnessBound}  yields  
$$\limsup_n\E\int_0^T[|\varphi_n(X^n_t)|+\psi_n(X^n_t)]^2dt<\infty,$$ which in turn implies  \eqref{tightAC}.
  \end{remark}

\subsection{Dol\'eans exponential and ``goodness''}
For a c\`adl\`ag semi--martingale $X = (X_t, t\geq 0)$, consider the stochastic linear equation of Dol\'eans
\begin{equation}\label{dolean}
Y_t= 1+ \int_0^t Y_{r-}dX_r .
\end{equation}
The following proposition follows from Theorem 1 and Theorem 2 in \cite{LiSh}, page 122.
\begin{proposition}\label{expdol}
{\rm (1)} Equation \eqref{dolean} has a unique solution (up to indistinguishability) within the class of semi--martingales. This solution is denoted by $\mathcal{E}(X)$ and is called the Dol\'eans exponential of $X$. It has the following representation
\begin{equation*}
\mathcal{E}(X)_t= \exp\big\{X_t-X_0-\frac12\langle X^c\rangle_t\big\}\prod_{r\leq t}(1+\Delta X_r)e^{-\Delta X_r}.
\end{equation*}
{\rm (2)} If $U$ and $X$ are two semi--martingales, then
\begin{equation*}
\mathcal{E}(U)_t\, \mathcal{E}(X)_t= \mathcal{E}(U+X+[U,X])_t.
\end{equation*}
{\rm (3)}
If $X$ is a local martingale, then $\mathcal{E}(X)$ is a nonnegative local martingale and a super--martingale.
\end{proposition}
\vspace{1cm}
For $\delta >0$ we define $h_\delta :\mathbb{R}_{+}\rightarrow \mathbb{R}_{+}$ by $h_\delta (r)= (1-\delta /r)^{+}$. For $x \in D([0,\infty ))$, we define $x^\delta \in D([0,\infty ))$ by 
$$x^\delta_t := x_t-\sum _{0<s\leq t}h_\delta (\mid \Delta x_s\mid )\Delta x_s.$$

\begin{definition}\label{defUCV}
{\rm(1)} Let $G, G^n$ in $\mathbb{D}, \{G^n, n\geq 1\}$ be a sequence of semimartingales adapted to a given filtration $(\mathcal{F}_t)$  and assume $G^n\Rightarrow G$ as $n\to\infty$. The sequence 
$(G^n)$ is called {\em good} if for any sequence $\{I^n, n\geq 1\}$ of $(\mathcal{F}_t)$--progressively measurable processes in $\mathbb{D}$ such that $(I^n, G^n)\Rightarrow (I, G)$ as $n\to\infty$, then $G$ is a semi--martingale for a filtration with respect to which $I$ is adapted, and 
$(I^n, G^n,\int I_{s-}^n dG_s^n) \Rightarrow (I,G,\int I_{s-}d G_s)$ as $n\to\infty$.

\noindent
{\rm(2)} A sequence of semi--martingales $\{G^n\}_{n\geq 1}$ is said to have {\em uniformly controlled variations} if there exists $\delta >0$, and for each $\alpha >0, n\geq 1$, there exists a semi--martingale decomposition $G^{n,\delta }= M^{n,\delta }+ A^{n,\delta }$ and a stopping time $T^{n,\alpha }$ such that $\mathbb{P}(\{T^{n,\alpha }\leq \alpha \})\leq \frac1{\alpha }$ and furthermore
\begin{equation}\label{UCV}
\sup _n \mathbb{E} \Big\{[M^{n,\delta },M^{n,\delta }]_{t \land T^{n,\alpha }}+ \int _0^{t \land T^{n,\alpha }} \mid dA^{n,\delta }\mid \Big\}<\infty .
\end{equation}
\end{definition}
It follows from pages 32 ff. in \cite{KP}
\begin{proposition}\label{th-good}
Let $G, G^n$ in $\mathbb{D}, \{G^n, n\geq 1\}$ be a sequence of semi--martingales and assume $G^n\Rightarrow G$.
\begin{description}
\item[\rm{(1)}]The sequence $\{G^n\}$ is good if and only if it has uniformly controlled variations.
\item[\rm{(2)}]If $\{G^n\}$ is good, then $(G^n, \mathcal{E}(G^n))\Rightarrow (G, \mathcal{E}(G)).$
\item[\rm{(3)}]Suppose $(I^n,G^n)\Rightarrow (I,G)$, and $\{G^n\}$ is good. Then $J^n= \int I_{s-}^n dG_s^n$, $n=1,2,\ldots$, is also a good sequence of semi--martingales. Moreover under the same conditions, $(I^n,G^n,\mathcal{E}(G^n))\Rightarrow
(I,G, \mathcal{E}(G))$.
\end{description}
\end{proposition}
\subsection{Two Girsanov theorems}
We state two multivariate versions of the Girsanov theorem, one for the Brownian and one for the point process case. The second one combines Theorems T2 and T3 from \cite{PB}, pages 165--166.
\begin{proposition}\label{Gir-BM}
Let $\{(B_s^{(1)},\ldots, B_s^{(d)}),\ s\ge0\}$  be a $d$--dimensional standard Brownian motion defined on the filtered  probability space $(\Omega ,\mathcal{F},\mathbb{P})$. Moreover,  let $ {\phi} = (\phi _1,...,\phi _d)$ be an $\mathcal{F}$-progressively measurable process with $\int_0^s \phi _i(r)^2 dr < \infty$ for all $1\le i\le d$ and $s\ge 0$. Let $X^{(i)}_s := \int_0^s \phi_i (r) \, dB^{(i)}_r$ and put $Y:= \mathcal E(X^{(1)} + \cdots + X^{(d)})$, or in other words
\begin{align*}
Y_s=  \exp \Big\{ \int_0^s \langle \phi (r), dB_r \rangle- \frac12 \int_0^s \mid \phi (r)\mid ^2 dr \Big\}.
\end{align*}
If $\mathbb{E}[Y_s]= 1$, $s\ge 0$, then    $\tilde{B}_s:=  B_s- \int_0^s \phi (r) dr$, $s\ge 0$, is a $d$-dimensional standard Brownian motion under  the probability measure $\tilde {\mathbb{P}}$ defined by 
$d{\tilde {\mathbb{P}}\mid_{\mathcal F_s}} / d{ {\mathbb{P}}\mid_{\mathcal F_s}} = Y_s$, $s\ge 0$.
\end{proposition}

\begin{proposition}\label{Gir-PP}
Let $\{(Q_s^{(1)},...,Q_s^{(d)}), s\ge 0\}$ be a $d$-variate point process adapted to some filtration $\mathcal{F}$, and let $\{\lambda _s^{(i)}, {s\ge 0}\}$  be the predictable $(\mathbb{P},\mathcal{F})$-intensity of $Q^{(i)}, 1\leq i\leq d$. Assume that none of the $Q^{(i)}$, $Q^{(j)}$, $i\neq j$, jump simultaneously. 
Let $\{\mu _r^{(i)}, r\ge 0\},  1\leq i\leq d$, be nonnegative $\mathcal{F}$-predictable processes such that for all $s\geq 0$ and all $1\leq i\leq d$
\begin{equation*}
\int _0^s \mu _r^{(i)} \lambda _r^{(i)} dr <\infty \qquad \mathbb{P} \mbox{ -a.s.}
\end{equation*}
For $i=1,\ldots,d$ and $s\ge 0$ define
$$X_s^{(i)} := \int_0^s (\mu_r^{(i)}-1)dM_r^{(i)}, \quad Y^{(i)}:= \mathcal E(X^{(i)}), \quad 
Y= \mathcal E(X^{(1)}+\ldots + X^{(d)}).
$$
Then, with $\{T_k^i, k=1,2\ldots\}$ denoting the jump times of $Q^{(i)}$,
$$ Y^{(i)}_s=  \Big( \prod _{k\geq 1: T_k^i\leq s} \mu^{(i)} _{T_k^i} \Big) \exp \Big\{ \int _0^s (1-\mu^{(i)} _r)\lambda _r^{(i)} dr \Big\} \quad  \mbox{ and } \quad
Y_s= \prod _{j=1}^d Y^{(j)}_s, \quad s\ge 0.$$
If $\mathbb{E}[Y_s]= 1$, $s \ge 0$, then, for each $1\leq i\leq d$, the process $Q^{(i)}$ has the $(\tilde {\mathbb{P}}, \mathcal{F})$-intensity $\tilde {\lambda }_r^{(i)}= \mu _r^{(i)} \lambda _r^{(i)} $, $r \ge 0$,
where the probability measure $\tilde {\mathbb{P}}$ is defined by 
$d{\tilde {\mathbb{P}}\mid_{\mathcal F_s}} / d{ {\mathbb{P}}\mid_{\mathcal F_s}} = Y_s$, $s\ge 0$.
\end{proposition}
\paragraph{Acknowledgements.} We thank  J.F. Delmas and Ed Perkins for stimulating discussions, J.F. Le Gall for calling our attention to the reference \cite{NRW}, and  a referee whose detailed report helped us to significantly improve the exposition.
{}
 

\begin{thebibliography}{99}
\addcontentsline{toc}{section}{References}
 \bibitem{BPS}  M. Ba, E. Pardoux, A. B. Sow,   Binary trees, exploration processes, and an
extended Ray--Knight Theorem,  {\it J. Appl. Probab.} {\bf49}, 2012.
\bibitem{PBi} P. Billingsley, {\it Probability and measure}, 3rd. ed., John Wiley and Sons Inc., New York, 1995.
\bibitem{Bi} P. Billingsley, {\it Convergence of Probability Measures}, 2d ed., John Wiley and Sons Inc., New York, 1999.
\bibitem {PB} P. Br\'emaud, {\it Point processes and queues: martingale dynamics}, Springer-Verlag New York, 1981.
\bibitem{EK} S. Ethier, Th. Kurtz, {\it  Markov processes: characterization and convergence}, John Wiley and Sons Inc., New York, 1986.
\bibitem{Fr} A. Friedman, {\it  Stochastic differential equations and applications, vol 1}, Academic Press, 1975.
\bibitem{JoMe} A. Joffe, M. M\'etivier, Weak convergence of sequences of semi--martingales with applications to multitype branching processes, {\it Adv. Appl. Prob.} {\bf18} (1986), 20--65.
\bibitem{KS} I. Karatzas, S. Shreve, {\it Brownian motion and stochastic calculus}, Graduate Texts in Mathematics {\bf113}, Springer, 1988.
\bibitem {KP} Th.  Kurtz and Ph. Protter, Weak convergence of stochastic integrals and differential equations, {\it 
Probabilistic models for nonlinear partial differential equations, Lecture Notes in Math} {\bf 1627}, 1--41, 1996.
    \bibitem{AL} A. Lambert, The branching process with logistic growth, {\it
    Ann. Appl. Probab.} {\bf15} (2005), 1506--1535.
\bibitem{LG4} J-F. Le Gall, It\^o's excursion theory and random trees,  \textit{Stochastic Process. Appl.} \textbf{120} (2010), 721--749.
    \bibitem {LiSh} R. S. Liptser and A. N. Shiryayev, {\it Theory of martingales}, Kluwer Academic Publishers, 1989.
      \bibitem {SM} S. M\'el\'eard, Quasi-stationary distributions for population processes, {\it Lecture at CIMPA school, St Louis, S\'en\'egal}, 2010, http://www.cmi.univ-mrs.fr/$\sim$pardoux/Ecole\_CIMPA/CoursSMeleard.pdf 
    \bibitem{NRW} J. R. Norris, L. C. G. Rogers, D. Williams, Self-avoiding random walk: a Brownian motion model with local time drift, {\it Probab. Th. Rel. Fields} {\bf 74} (1987), 271-287.
   \bibitem{PW} E. Pardoux and A. Wakolbinger, From exploration paths to mass excursions - variations on a theme of Ray and Knight,  in {\it Surveys in Stochastic Processes}, Proceedings of the 33rd SPA Conference in Berlin, 2009, J. Blath, P. Imkeller, S. Roelly (eds.), pp. 87--106, EMS 2011.
 \bibitem{PW2}  E. Pardoux and A. Wakolbinger, From Brownian motion with a local time drift to Feller's branching diffusion with logistic growth,  {\it Elec. Comm. Probab.}, to appear.
 \bibitem{Pe} E. Perkins, Weak invariance principles for local time, {\it Z. Wahrscheinlichkeitstheorie verw. Gebiete} {\bf 60} (1982), 437--451.
\bibitem{Pi} J. Pitman, The distribution of local times of a Brownian bridge, 
S\'eminaire de probabilit\'es (Strasbourg) {\bf33}, Lecture Notes in Math. {\bf1709}, pp. 388--394, 1999.
 \bibitem{ReYo} D. Revuz and M. Yor, {\it Continuous martingales and Brownian motion}, 3rd ed., Spinger Verlag, New York   1999. 
 \bibitem{DWS} D. W. Stroock, {\it Probability theory: an analytic view}, Cambridge University Press, 1993.
\end{thebibliography}
 \end{document}